\documentclass[11pt,a4paper]{article}

\usepackage{authblk}
\usepackage{color}
\usepackage{natbib}
\usepackage[utf8]{inputenc} 
\usepackage[T1]{fontenc}    
\usepackage{textcomp}
\usepackage{hyperref}       
\usepackage{url}            
\usepackage{booktabs}       
\usepackage{amsfonts}       
\usepackage{nicefrac}       
\usepackage{microtype}      
\usepackage{amsmath,amsfonts,amssymb}
\usepackage{algorithm}
\usepackage{algorithmic}
\usepackage{mathtools}
\usepackage{graphicx}
\usepackage{amsthm}
\usepackage{caption}
\usepackage{subcaption}
\usepackage{subfiles}
\newcommand{\norm}[1]{\left\lVert#1\right\rVert}
\newcommand\inner[2]{\langle #1, #2 \rangle}
\DeclareMathOperator*{\argmin}{arg\,min}
\newtheorem{theorem}{Theorem}
\newtheorem{corollary}{Corollary}
\newtheorem{prop}{Proposition}

\newcommand{\x}{\mathbf{x}}
\newcommand{\y}{\mathbf{y}}
\newcommand{\z}{\mathbf{z}}
\newcommand{\s}{\mathbf{s}}
\renewcommand{\L}{\mathbf{L}}

\title{\bf Communication-Efficient Projection-Free Algorithm for Distributed Optimization}

\author[1]{Yan Li}
\author[2]{Chao Qu}
\author[1]{Huan Xu}
\affil[1]{H. Milton Stewart School of Industrial and Systems Engineering, Georgia Institute of Technology}
\affil[2]{Department of Electrical Engineering, Technion }

\begin{document}

\maketitle

\begin{abstract} Distributed optimization has gained a surge of interest in
recent years. In this paper we propose a distributed projection free algorithm
named Distributed Conditional Gradient Sliding(DCGS). Compared to the state-of-the-art
distributed Frank-Wolfe algorithm, our algorithm attains the same communication
complexity under much more realistic assumptions. In contrast to the consensus based algorithm, DCGS is based on the primal-dual algorithm, yielding a modular analysis that
can be exploited to improve linear oracle complexity whenever centralized
Frank-Wolfe can be improved. We demonstrate this advantage and show that the
linear oracle complexity can be reduced to almost the same order of magnitude
as the communication complexity, when the feasible set is polyhedral. Finally we
present experimental results on Lasso and matrix completion, demonstrating significant
performance improvement compared to the existing distributed Frank-Wolfe algorithm. \end{abstract}

\section{Introduction}
Decentralized optimization methods have been widely used in the machine learning community. Compared to  centralized optimization methods, they
enjoy several advantages, including aggregating the computing power of distributed
machines, robustness to dynamic network topologies, and preserving data privacy \citep{yin_decentralized_grad_descent}. These
advantages make them an attractive option when data are collected by
distributed agents, and either communicating the data to a central processing
agent is computational prohibitive, or due to privacy concern each agent needs to
keep the local data privately.

After the seminal work in  \citet{tsitsiklis_seminal_2,tsitsiklis_seminal_1},
there have been fruitful development on  distributed optimization. For non-smooth
objective functions, consensus-based subgradient methods have been analyzed in
\citet{nedic_multi_agent,nedic_stochastic_subgrad_projection,nedic_aynchronous}.
Dual averaging method proposed in \citep{duchi_dual_averaging}
further explains the effect of the network topology on the convergence. Usually
distributed subgradient methods converge at the rate of
$\mathcal{O}(\frac{\log(T)}{\sqrt{T}})$. For smooth
(and possibly strongly convex) objective, \citet{yin_extra,yin_proximal_extra}
propose an \emph{exact} first-order algorithm EXTRA and its proximal variant,
which attain an improved rate of $\mathcal{O}(\frac{1}{T})$
for general smooth objective, and a linear rate of $\mathcal{O}(c^T)$ for smooth and strongly
convex objective with $c < 1$ (see also
\citet{nedic_geometric_convergence,yin_decentralized_grad_descent}).
Asynchronous decentralized (sub)gradient descent algorithms have also been proposed and analyzed
\citet{nedic_aynchronous,tsitsiklis_seminal_1}. Another class of mainstream distributed algorithms are based
on the dual method, which include the classic idea of dual decomposition
\citep{murray_dual_decomposition}, and the celebrated alternating direction
method of multipliers (ADMM)
\citep{boyd_admm,ozdaglar_distributed_admm,ozdaglar_asynchronous_dadmm}. ADMM
attains $\mathcal{O}(\frac{1}{T})$ convergence rate for the smooth problem and
$\mathcal{O}(c^T)$ for smooth and strongly convex problem, but such results
rely on strong assumptions such as having no constraints
\citep{yin_linear_convergence_admm} or the local subproblem of each agent being easy
to solve \citep{ozdaglar_asynchronous_dadmm}.

All the aforementioned methods can be categorized as \emph{projection-based}
methods, as they all require to take projection back to the feasible set of the constraints at each iteration. Though
commonly assumed to be easy,  in numerous applications such
projection indeed can either be computational expensive  (projection onto the trace norm ball, base
polytopes \citep{Fujishige_asubmodular}) or even
intractable \citep{collins_intractable}. Frank-Wolfe (FW) algorithm arises as a
natural alternative in these scenarios. Unlike projection-based methods, FW
assumes a linear oracle (LO) that solves a linear optimization problem over the
feasible set which can be significantly easier than the projection. We
refer to algorithms that avoid projection as \emph{projection-free}
algorithms. FW algorithm has been revisited in recent years for its projection-free property
\citep{jaggi_revisiting_fw,jaggi_phd_thesis,clarkson_coreset,bach_herding_and_fw} and numerous improvements have been made. These include
regularized FW \citep{nemirovski_regularzied_fw}, linearly convergent variants under additional assumptions
\citep{hazan_linear_convergence,julien_linear_convergence},  and stochastic variants
\citep{hazan_proj_free_online_learning,hazan_vr_projection_free}.

\textbf{Related Work.} Despite the progress on centralized FW algorithm, results on distributed FW
algorithm are surprisingly limited. Specialized versions of decentralized FW
algorithm have been proposed. \citet{wang_pdfw} propose a distributed
block-coordinate FW algorithm for block-separable feasible sets
\citep{julien_bcfw}. \citet{balcan_dfw_sparse_learning} consider a Lasso-type
distributed learning problem. They neither assume nor exploit the fact that
the global objective is a natural summation of each agent's local objective,  and their communication scheme is also different from the usual ``within
neighborhood'' communication scheme and could be significantly more complicated.
\citet{wai_dfw_comm_efficient} consider a distributed FW algorithm also for
the Lasso-type problem that leverages the sparsity of iterates to improve
communication overhead. To the best of our knowledge, the most recent
distributed FW (DFW) algorithm on general smooth convex problems is by
\citet{wai_dfw_convex_nonconvex}, wherein the DFW convergence rate is
$\mathcal{O}(\frac{1}{T})$ for smooth objectives; and
$\mathcal{O}(\frac{1}{T^2})$ for smooth and strongly convex objective {\em under the assumption that  the minimizer lies in the  interior of constraint set}. This assumption is almost unrealistic for two reasons: it implies the problem is essentially unconstrained, which usually fails to impose structural properties (such as sparsity, low-rankness) to the solution; and for a unconstrained problem the vanilla distributed gradient descent algorithm suffices to solve the problem efficiently \citep{yin_decentralized_grad_descent}. Whether such
restrictive assumption could be removed  while retaining the same complexity remains an open question.

We should note that all the previously discussed methods share the
same communication complexity and projection/LO complexity to obtain an
$\epsilon$-optimal solution, regardless of being projection-based or
projection-free. In practice, however, the time consumed by a single communication
and a LO/projection often differ by \emph{orders of magnitude}, which could incur
significant latency. Modern CPUs perform IO at over 10 GB/s yet
communication over TCP/IP is about 10 MB/s, this gap is even more significant
when LO oracle is already cheap. Consider the matrix completion problem, in Section \ref{experiments} we will show that for an iteration consisting one round of communication and one LO, communication would take up over $97 \%$ time. This implies in DFW  the actual running time   would be largely consumed by communication.
To alleviate the problem of latency in communication expensive applications, whether it is possible to trade for a better
communication complexity with a moderately worse LO complexity becomes another open question.

\textbf{Contributions.} In this paper we answer the above mentioned questions with an (almost)
affirmative answer. Our contributions are the following:
\begin{itemize}
  \item We propose a new distributed projection-free algorithm named Decentralized
  Conditional Gradient Sliding (DCGS), and show that it attains
  $\mathcal{O}(\frac{1}{\epsilon})$ communication complexity and
  $\mathcal{O}(\frac{1}{\epsilon^2})$ LO complexity for smooth convex objectives.

 \item Without assuming the minimizer being in the {\em interior} of the feasible set, we
 show that DCGS attains $\mathcal{O}(\frac{1}{\sqrt{\epsilon}})$ communication
 complexity and $\mathcal{O}(\frac{1}{\epsilon^2})$ LO complexity for
 smooth and strongly convex objectives.

 \item Our algorithm builds upon a distributed version of primal-dual algorithm
 and is hence modular. As a consequence, improvement on centralized FW can  be easily
 exploited by DCGS. We demonstrate this advantage when the feasible
 set is polyhedral, for smooth and convex objective the LO complexity can be reduced to
 $\tilde{\mathcal{O}} (\frac{1}{\epsilon})$\footnote{Throughout this paper we use
 $\tilde{\mathcal{O}} () $ to hide any additional logrithmic factor}, while for smooth and strongly convex
 objective the LO complexity can be further reduced to $\tilde{\mathcal{O}}
 (\frac{1}{\sqrt{\epsilon}})$, which matches the result of
 \citep{wai_dfw_convex_nonconvex}, but without the restrictive assumption on the minimizer.

 \end{itemize}

\section{Problem formulation}
 We consider an undirected graph $G = (V, E)$, where $V = \{1, \ldots, m \}$ denotes the vertex set and $E \subset V \times V$ denotes the edge set. Each node $i \in V$ is associated with an agent indexed also by $i$, and has its local objective $f_i(x): \mathbb{R}^d \rightarrow \mathbb{R}$. We define the neighborhood of agent $i$ to be $N(i) = \{j \in V: (i,j) \in E \}$. Each agent $i$ can only communicate information with its neighbors. Naturally, $G$ is assumed to be connected or otherwise distributed optimization is impossible. Our objective is to minimize the summation of the local objectives, subject to the constraint that $x$ belongs to a closed compact convex set $X$, that is:
\begin{align}\label{distributed_formulation}
 \underset{x \in X}{\mbox{Minimize:}} \quad f(x):= \sum_{i=1}^m f_i(x).
\end{align}
We assume each function $f_i$ is $u$ (possibly $0$)-strongly convex and $l$-smooth,
i.e., $ \frac{u}{2} \norm{y-x}^2 \leqslant f_i(y) - f_i(x) - \nabla f_i(x)^{T} (y-x) \leqslant \frac{l}{2} \norm{y-x}^2 $. Our algorithm could also be easily adapted to the setting where $f_i$ has different smoothness and strong convexity. We present here only the homogeneous case for simplicity.

The distributed formulation (\ref{distributed_formulation}) can be reformulated as the following linearly constrained optimization problem. Consider each agent $i$ keeps its local copy of decision variable $x_i$, we can impose a linear constraint on $\mathbf{x} = (x_1, \ldots, x_m)$ so that $x_i = x_j$ for all $(i,j) \in E$. Define the graph Laplacian $L \in \mathbb{R}^{m \times m}$ to be:
\begin{align}
 L_{ij} = \begin{cases}
  |N(i)| & \text{if }  i=j, \\
  -1 & \text{if }  i \neq j \text{ and } (i,j) \in E, \\
  0 &  \text{otherwise.}
\end{cases}
\end{align}
Then (\ref{distributed_formulation}) could be reformulated as:
\begin{align} \label{linear_constraint_formulation}
  \begin{split}
 \min_{\mathbf{x} \in X^m}\quad F(\mathbf{x})  & = \sum_{i = 1}^m f_i(x_i) \\
 s.t. \quad \mathbf{Lx}  &=  \mathbf{0},
\end{split}
\end{align}
Let $\otimes$ denotes the Kronecker product, here $\mathbf{L} = L \otimes I_d$ and $X^m = \{(x_1,\ldots, x_m): x_i \in X \}$. Since $G$ is assumed to be connected,   (\ref{distributed_formulation}) and (\ref{linear_constraint_formulation}) are equivalent.

We can further reformulate the linearly constrained problem as a bilinear saddle point problem. Observe that (\ref{linear_constraint_formulation}) is equivalent to:
\begin{align} \label{bilinear_saddle_formulation}
  \min_{\mathbf{x} \in X^m} \max_{\mathbf{y} \in \mathbb{R}^{md}} F(\mathbf{x}) + \inner{\mathbf{Lx}}{\mathbf{y}}.
\end{align}
The bilinear saddle point problem (\ref{bilinear_saddle_formulation}) is well suited for the primal-dual algorithm proposed in \citep{chambolle_primal_dual}. We present the orginal primal-dual algorithm applied to our problem in Algorithm~\ref{alg:raw_primal-dual}.

\begin{algorithm}
    \caption{Primal-dual algorithm \citep{chambolle_primal_dual}}
    \label{alg:raw_primal-dual}
    \begin{algorithmic}[1]
        \STATE{Let $\mathbf{x}^{0} = \mathbf{x}^{-1} \in X^m$ and $\mathbf{y}^0 \in \mathbb{R}^{md}$ and $\{\alpha_k\}, \{\tau_k\}, \{\eta_k\}$ be given.}
        \FOR{$k=1, \ldots, N$ do}
        \STATE{ $\tilde{\mathbf{x}}^k = \alpha_k (\mathbf{x}^{k-1} - \mathbf{x}^{k-2}) + \mathbf{x}^{k-1}$}
        \STATE{ $ \mathbf{y}^k = \argmin_{\mathbf{y} \in \mathbb{R}^{md}} \inner{-\mathbf{L} \tilde{\mathbf{x}}^k}{\mathbf{y}} + \frac{\tau_k}{2} \norm{\mathbf{y} - \mathbf{y}^{k-1}}^2$}
        \STATE{ $ \mathbf{x}^k = \argmin_{\mathbf{x} \in x^m} \inner{\mathbf{L y^k}}{\mathbf{x}} + F(\mathbf{x}) + \frac{\eta_k}{2} \norm{\mathbf{x}^{k-1} - \mathbf{x}}^2 $}
        \ENDFOR
    \end{algorithmic}

\end{algorithm}

\citet{lan_comm_efficient} observe that since $F(\mathbf{x})$ is a summation splitted across agents, all the updates in the primal-dual algorithm can be performed in a distributed way. They propose a distributed primal-dual algorithm and show that to find an $\epsilon$-optimal solution, one needs $\mathcal{O}(\frac{1}{\epsilon})$ rounds of communication for a non-smooth convex objective and $\mathcal{O}(\frac{1}{\sqrt{\epsilon}})$ for a non-smooth strongly convex objective, which improves upon the previous results. However, their algorithm still lies in the category of projection-based algorithms and they consider non-smooth problem, which is different from our setting.


\section{Algorithm}
In this section we present in Algorithm \ref{alg:pdcgs} the Decentralized Conditional Gradient Sliding for a general convex feasible set equipped with a linear oracle.
At a high-level, DCGS is closely related to Conditional Gradient Sliding (CGS) algorithm proposed in \citet{lan_cgs}. However CGS considers only the primal problem, here we consider a primal-dual problem due to performing distributed optimization. As such, the analysis is significantly more involved.

\begin{algorithm}
    \caption{Distributed Conditional Gradient Sliding (DCGS)}
    \label{alg:dcgs}
    \begin{algorithmic}[1]
        \STATE{Let $\mathbf{x}^{0} = \mathbf{x}^{-1} \in X^m$ and $\mathbf{y}^0 \in \mathbb{R}^{md}$ and $\{\alpha_k\}, \{\tau_k\}, \{\eta_k\}$ be given.}
        \FOR{$k=1, \ldots, N$ do}
        \STATE{Update for all agents as the following:}
        \STATE{ $\tilde{x}_i^k = \alpha_k (x_i^{k-1} - x_i^{k-2}) + x_i^{k-1} $}
        \STATE{ $ v_i^k = \sum_{j \in N(i) \cup \{i\}} L_{ij} \tilde{x}_j^k $}
        \STATE{ $ y_i^k = y_i^{k-1} + \frac{1}{\tau_k} v_i^k $}
        \STATE{ $ w_i^k = \sum_{j \in N(i) \cup \{i\} } L_{ij} y_j^k$}
        \STATE{ $x_i^k = CG(f_i, x_i^{k-1}, w_i^k, \eta_k, e_i^k)$}\label{lo_subproblem}
        \ENDFOR
        \RETURN $\overline{\mathbf{x}}^N = (\sum_{k=1}^N \theta_k)^{-1} \sum_{k=1}^N \theta_k \mathbf{x}^k$
        \STATE{}
        \STATE{\textbf{Procedure:} $x^+ = CG(f, x, w, \eta, e)$}\label{cg_procedure}
        \STATE{Let $x^0 \in X$}
        \WHILE{$t=0,\ldots$}
        \STATE{Let $s^t = \argmin_{s \in X} \inner{\nabla f(x^t) + w + \eta(x^t - x)}{s}$}
        \IF{ $\inner{\nabla f(x^t) + w + \eta(x^t - x)}{x^t - s^t} \leqslant e$ }
        \RETURN $x^t$
        \ENDIF
        \STATE{ $x^{t+1} = (1 - \gamma_t) x^t + \gamma_t s^t$, where $\gamma_t  = \frac{2}{t+2}$ or computed by line search}
        \ENDWHILE
        \STATE{\textbf{End procedure}}
    \end{algorithmic}
\end{algorithm}
The most important step of DCGS algorithm is in Line \ref{lo_subproblem}, where we update decision variable $x_i^k$ by calling the CG procedure defined in Line \ref{cg_procedure}. If we define $\phi_i^k (x_i) = \inner{w_i^k}{x_i} + f_i(x_i) + \frac{\eta_k}{2} \norm{x_i - x_i^{k-1}}^2$, then the CG procedure could be seemed as the FW algorithm applied to $\min_{x_i \in X} \phi_i^k(x_i)$, with termination criterion $ \inner{ \nabla \phi_i^k(x_i^k)}{x_i^k - x_i^{k-1}} \leqslant e_i^k$, where the left hand side is often called the Wolfe gap. Below we make a few remarks on the communication mechanism, the main technical challenges and the modularity of our algorithm.

\textbf{Communication Mechanism.}
At each outer iteration $k$, each local agent first computes $\tilde{x}_i^k$ based on
extrapolation of two previous primal iterates, and broadcase $\tilde{x}_i^k$ to
all of its neighbors $j \in N(i)$. After one round of broadcasting, each agent
uses $\tilde{x}_j^k$ received from its neighbors and perform dual variable update
$y_i^k$, then broadcast the updated dual variable to all of its neighbors $j \in
N(i)$. After second round of broadcasting, each agent uses $y_j^k$ received from
its neighbors and call the CG procedure to update primal variable $x_i^{k+1}$. Each iteration incurs two rounds
of communication within the network, hence the overall communication complexity is
the same as the outer iteration complexity.

\textbf{Trade-off between Communication and LO.}
If we set $e_i^k=0$ in Line \ref{lo_subproblem}, we are solving problem $\min_{x_i \in X} \phi_i^k(x_i)$ exactly. The outer iteration of DCGS then reduces to the primal-dual algorithm applied to our problem (\ref{bilinear_saddle_formulation}), implemented in a distributed fashion. By well-known results of the primal-dual algorithm \citep{chambolle_primal_dual,lan_comm_efficient}, to get an $\epsilon$-optimal solution one needs $\mathcal{O}(\frac{1}{\epsilon})$ iterations for a convex objective and $\mathcal{O}(\frac{1}{\sqrt{\epsilon}})$ iterations for a strongly convex objective.
From our previous discussion, this yields the communication complexity of DCGS to be $ \mathcal{O}(\frac{1}{\epsilon})$ and $\mathcal{O}(\frac{1}{\sqrt{\epsilon}})$ respectively. However in this extreme case LO calls in the CG procedure would be prohibitively large.
Consequently, we need to carefully choose $e_i^k$ to ensure that the subproblem $\min_{x_i \in X} \phi_i^k(x_i)$ is solved in a controlled way: the convergence of outer iteration should be approximately at the same speed as the case when the subproblem is solved exactly, but meanwhile we need to keep LO complexity in the CG procedure to remain relatively small.

\textbf{Modularity.} The CG procedure in DCGS algorithm is where all the calls to LO take place.  We believe there are not much room for improvement in terms of the communication complexity, as our complexity in the general case matches that of the DFW algorithm under additional (very strong) assumption that the optimal solution is in the interior of the feasible set. The room for improvement then lies in possibly reducing LO complexity. If we treat the CG procedure as a module in the DCGS algorithm, can we replace it with a module that runs much faster for specific objectives or feasible sets, and obtain a better DCGS variant? The answer is affirmative. As an example we will show that significant improvement on LO complexity could be made when the feasible set is polyhedral.

\section{Theoretical Results}\label{result}
\subsection{General Feasible Set}
In this section we set suitable parameters to DCGS for convex and strongly convex objectives. We will present its convergence results, communication and LO complexity. We also present a detailed comparison with results of DFW in \citep{wai_dfw_convex_nonconvex}.

\begin{theorem}[Convergence for Smooth and Convex Objectives]\label{thrm:smooth}
  \hfill \\
  Set $ \theta_k = \alpha_k = 1, \eta_k = 2\norm{L}, \tau_k = \norm{L}, e_i^k = \frac{\norm{L} \max(\norm{\mathbf{x}^0 - \mathbf{x}^*}^2, \norm{\mathbf{y}^0}^2)}{mN}$ in Algorithm \ref{alg:dcgs}, where $\norm{L}$ denotes the spectral norm of Laplacian matrix $L$.
  Assume each $f_i$ is $l$-smooth, we have:
  \begin{align} \label{smooth_convergence}
    F(\overline{\mathbf{x}}_N) - F(\mathbf{x}^*) \leqslant \frac{ \norm{L}}{N} \max(\norm{\mathbf{x}^0 - \mathbf{x}^{*}}^2, \norm{\mathbf{y}^0}^2).
  \end{align}
\end{theorem}
From (\ref{smooth_convergence}) it is straightforward to establish communication complexity. Note that it only depends on the number of agents and the network topology, and is \emph{independent} of the objective function $f$, which is a feature that DFW does not have.
\begin{corollary}[Complexity for Smooth and Convex Objectives]\label{complexity:smooth}
\hfill \\
Under the same conditions as in Theorem \ref{thrm:smooth}, to get a solution such that $F(\overline{\mathbf{x}}_N) - F(\mathbf{x}^*) \leqslant \epsilon$, the number of communications and LO for each agent are respectively bounded by:
\begin{align}\label{communication:smooth}
\mathcal{O} \left(\frac{\norm{L} \max(\norm{\mathbf{x}^0 - \mathbf{x}^{*}}^2, \norm{\mathbf{y}^0}^2)}{\epsilon} \right)
\end{align}
and
\begin{align}\label{lo:smooth}
 \mathcal{O} \left( \frac{ m \norm{L}(l + \norm{L})  \max(\norm{\mathbf{x}^0 - \mathbf{x}^{*}}^2, \norm{\mathbf{y}^0}^2)}{\epsilon^2}  \right)
\end{align}
\end{corollary}

\textbf{Detailed Comparison.} DFW in
\citep{wai_dfw_convex_nonconvex} has communication and LO complexity
both bounded by $\mathcal{O} \left(\frac{l m D^2 G}{\epsilon} \right)$, where
$D$ denotes the diameter of the feasible set and $G$ inversely relates to the spectral gap
of the weighted communication matrix. If we set $\mathbf{y}^0 = \mathbf{0}$, and observe that
$\frac{\norm{\mathbf{x}^0 - \mathbf{x}^*}}{D^2} = \mathcal{O}(\frac{1}{m})$, then it
can be seen that our algorithm is at least $\mathcal{O}
(\frac{lG}{\norm{L}})$ times better in terms of the communication complexity, and at
most $\mathcal{O} (\frac{ m(l+\norm{L})}{\epsilon} \cdot
\frac{\norm{L}}{lG})$ worse in LO complexity. Suppose in application where the
agents network is set beforehand so that $(m, \norm{L}, G)$ be treated as constants, as
objective becomes increasingly non-smooth, our algorithm outperforms DFW by factor of $\mathcal{O}(l)$ in communication with $\mathcal{O}(\frac{1}{\epsilon})$ worsened LO complexity. For tasks that communication is time consuming but  LO is much cheaper (e.g., matrix completion), such a trade-off can be significant, especially when we are not solving for a high precision solution.

\begin{theorem}[Convergence for Smooth and Strongly Convex Objectives]\label{thrm:sc} \hfill \\
Set $\alpha_k = \frac{k}{k+1}, \theta_k = k+1, \eta_k = \frac{ku}{2}, \tau_k = \frac{4\norm{L}^2}{(k+1)u}, e_i^k = \frac{\max(u \norm{\mathbf{x}^0 - \mathbf{x}^*}^2, \norm{L}^2 \norm{\mathbf{y}^0}^2/u )}{mNk} $ in Algorithm \ref{alg:dcgs}. Assume $f_i$ is $u$-strongly convex and $l$-smooth, we have:
\begin{align}\label{sc:convergence}
 F(\overline{\mathbf{x}}_N) - F(\mathbf{x}^*) \leqslant \frac{1}{N^2} \max\left(u \norm{\mathbf{x}^0 - \mathbf{x}^*}^2, \frac{\norm{L}^2 \norm{\mathbf{y}^0}^2}{u} \right).
\end{align}
\end{theorem}
Similarly to establishing (\ref{communication:smooth}) and (\ref{lo:smooth}), we bound the communication and LO complexity in the following corollary.

\begin{corollary}[Complexity for Smooth and Strongly Convex Objective] \hfill \label{complexity:sc}\\
Under the same conditions as in Theorem \ref{thrm:sc}, to get a solution such that $F(\overline{\mathbf{x}}_N) - F(\mathbf{x}^*) \leqslant \epsilon$, the number of communications and LO for each agent can be respectively bounded by:
\begin{align}\label{communication:sc}
 \mathcal{O} \left( \sqrt{\frac{\max(u \norm{\mathbf{x}^0 - \mathbf{x}^*}^2, \frac{\norm{L}^2 \norm{\mathbf{y}^0}^2}{u} )} {\epsilon}}  \right),
\end{align}
and
\begin{align}\label{lo:sc}
\mathcal{O} \left( \frac{m l  \max(u\norm{\mathbf{x}^0 - \mathbf{x}^*}^2, \frac{\norm{L}^2 \norm{\mathbf{y}^0}^2}{u} )}{\epsilon^2} \right).
\end{align}
\end{corollary}

\textbf{Detailed Comparison.}
DFW in \citep{wai_dfw_convex_nonconvex} has both communication and LO complexity bounded by $\mathcal{O}(\frac{1}{\sqrt{\epsilon}})$, but requires the minimizer to be bounded way from boundary of the feasible set.
Our complexity result does not rely on this unrealistic assumption that often fails, especially when the constraint should be active to impose structural assumption (e.g., sparsity) on the solution. It is then fair to compare with their result in the convex setting. Our result can be deemed as trading for a better communication complexity from $\mathcal{O}(\frac{1}{\epsilon})$ to $\mathcal{O}(\frac{1}{\sqrt{\epsilon}})$, with a moderately worse LO complexity from $\mathcal{O}(\frac{1}{\epsilon})$ to $\mathcal{O}(\frac{1}{\epsilon^2})$.

\subsection{Polyhedral Feasible Sets}
The trade-off between communication and LO however is almost unnecessary, when the feasible set is polyhedral. Specifically, DCGS achieves the same communication and LO complexity (with additionally logarithmic factor), regardless of $f_i$ being strongly convex or not. This improvement is a direct result  from the modularity of DCGS: we replace the CG procedure in DCGS with a faster one adapted from \citep{julien_linear_convergence}. We present the modified DCGS for a polyhedral feasible set in Algorithm \ref{alg:pdcgs}.
\begin{algorithm}
    \caption{Distributed Conditional Gradient Sliding (DCGS) over polyhedral}
    \label{alg:pdcgs}
    \begin{algorithmic}[1]
      \STATE{ \dots as in Algorithm~\ref{alg:dcgs}, except replace Line \ref{lo_subproblem} by:}
      \STATE{$ x_i^k = PCG( f_i, x_i^{k-1}, w_i^k, \eta_k, e_i^k )$ }
      \STATE{}
      \STATE{\textbf{Procedure:} $x^+ = PCG(f,x,w,\eta,e)$}
      \STATE{Let $x^0 \in X, S^0 = \{x^0\}, \alpha_{x^0}^0 = 1 $ and $\alpha_{\mu}^0 = 0$ for all $\mu \in X\setminus \{x^0\}$}
      \WHILE{ $t=0, \ldots$}
      \STATE{ Compute $s^t = \argmin_{s \in X} \inner{\nabla f(x^t) + w + \eta(x^t - x)}{s}$}
      \STATE{ Compute $v^t = \argmin_{s \in S^t} - \inner{ \nabla f(x^t) + w + \eta(x^t - x)}{s}$}
      \IF{ $\inner{\nabla f(x^t) + w + \eta(x^t - x)}{x^t - s^t}  \leqslant e$}
      \RETURN $x^t$
      \ENDIF
      \STATE{Let $d^t = s^t -v^t$ and compute $\gamma_t = \argmin_{\gamma \in [0,\alpha_{v^t}^t]} f(x^t + \gamma d^t)$}
      \STATE{Update $x^{t+1} = x^t + \gamma_t d^t$}
      \STATE{Update $\alpha_{s^t}^{t+1} = \alpha_{s^t}^{t} + \gamma_t$, $\alpha_{v^t}^{t+1} = \alpha_{v^t}^{t} - \gamma_t $, and $\alpha_{\mu^t}^{t+1} = \alpha_{\mu^t}^{t}$ for all $\mu \in S^t \setminus \{s^t, v^t \}$ }
      \STATE{ Update $S^{t+1} = \{\mu \in X: \alpha_{\mu}^{t+1} >0 \}$}
      \ENDWHILE
    \end{algorithmic}
\end{algorithm}


\begin{corollary}[Smooth and Convex: polyhedral set] \label{complexity:smooth_poly}\hfill \\
Under the same conditions as in Theorem \ref{thrm:smooth}, suppose $X$ is a polyhedral set with pyramidal width $W$ and width $D$, for DCGS each agent has the same communication complexity as in (\ref{communication:smooth}), and has LO complexity:
\begin{align}\label{lo:smooth-poly}
 \tilde{\mathcal{O}} \left( \frac{D^2}{W^2} \frac{\log m(l+\norm{L}) \max( \norm{\mathbf{x}^0 - \mathbf{x}^*}^2, \norm{\mathbf{y}^0}^2)}{\epsilon} \right).
\end{align}
\end{corollary}

\textbf{Improvements.} Observe that the LO complexity is now at the same order of magnitude as in \citep{wai_dfw_convex_nonconvex} in terms of dependence on $\epsilon$. If we set $\mathbf{y}^0 = \mathbf{0}$, LO complexity of DCGS in this setting is at most $\tilde{\mathcal{O}} (\frac{D^2}{W^2} \frac{\log m (l+ \norm{L})}{lG})$ worse than in \citep{wai_dfw_convex_nonconvex}, which could even be better when the network is pooly connected (so that $G$ is large). Our complexity also depends on $\frac{D^2}{W^2}$ which could be interpretated as the condition number of the polyhedral constraint set $X$.

\begin{corollary}[Smooth and strongly convex: polyhedral set]\label{complexity:sc_poly}\hfill \\
  Under the same conditions as in Theorem \ref{thrm:sc}, suppose $X$ is a polyhedral set with pyramidal width $W$ and width $D$, for DCGS each agent has the same communication complexity as in (\ref{communication:sc}), and has LO complexity:
\begin{align}\label{lo:sc-poly}
 \tilde{\mathcal{O}} \left( \frac{D^2}{W^2} \frac{ l \log m }{u} \sqrt{\frac{\max( u \norm{\mathbf{x}^0 - \mathbf{x}^*}^2, \frac{\norm{L}^2 \norm{\mathbf{y}^0}^2}{u})}{\epsilon}}  \right).
\end{align}
\end{corollary}

\textbf{Improvements.} The LO complexity is now improved to $\tilde{\mathcal{O}} (\frac{1}{\sqrt{\epsilon}})$ which is of the same order of magnitude as in \citep{wai_dfw_convex_nonconvex} in terms dependence on $\epsilon$, but this result makes no assumption on the minimizer. If we choose $\mathbf{y}^0 = \mathbf{0}$, our LO complexity reduces to:
$
\tilde{\mathcal{O}} \left(\frac{D^2}{W^2} \frac{l \log m}{\sqrt{u}} \frac{\norm{\mathbf{x}^0 - \mathbf{x}^*}}{\sqrt{\epsilon}}  \right)
$
which has a clean interpretation in terms of its dependency on condition number of the objective and condition number of the polyhedral constraint.

\section{Experimental Results}\label{experiments}
In this section we present experiments comparing DCGS with the existing distributed FW algorithm in \citep{wai_dfw_convex_nonconvex} and demonstrate the superiority of our algorithm. For both of the following experiments, we set the number of agents $m=10$, and the associated network to be a $1$-connected cycle, i.e. each agent $i$ is connected to its previous one $i-1$ and the latter one $i+1$.

\begin{figure}
\begin{subfigure}{.26\textwidth}
  \includegraphics[width=\linewidth]{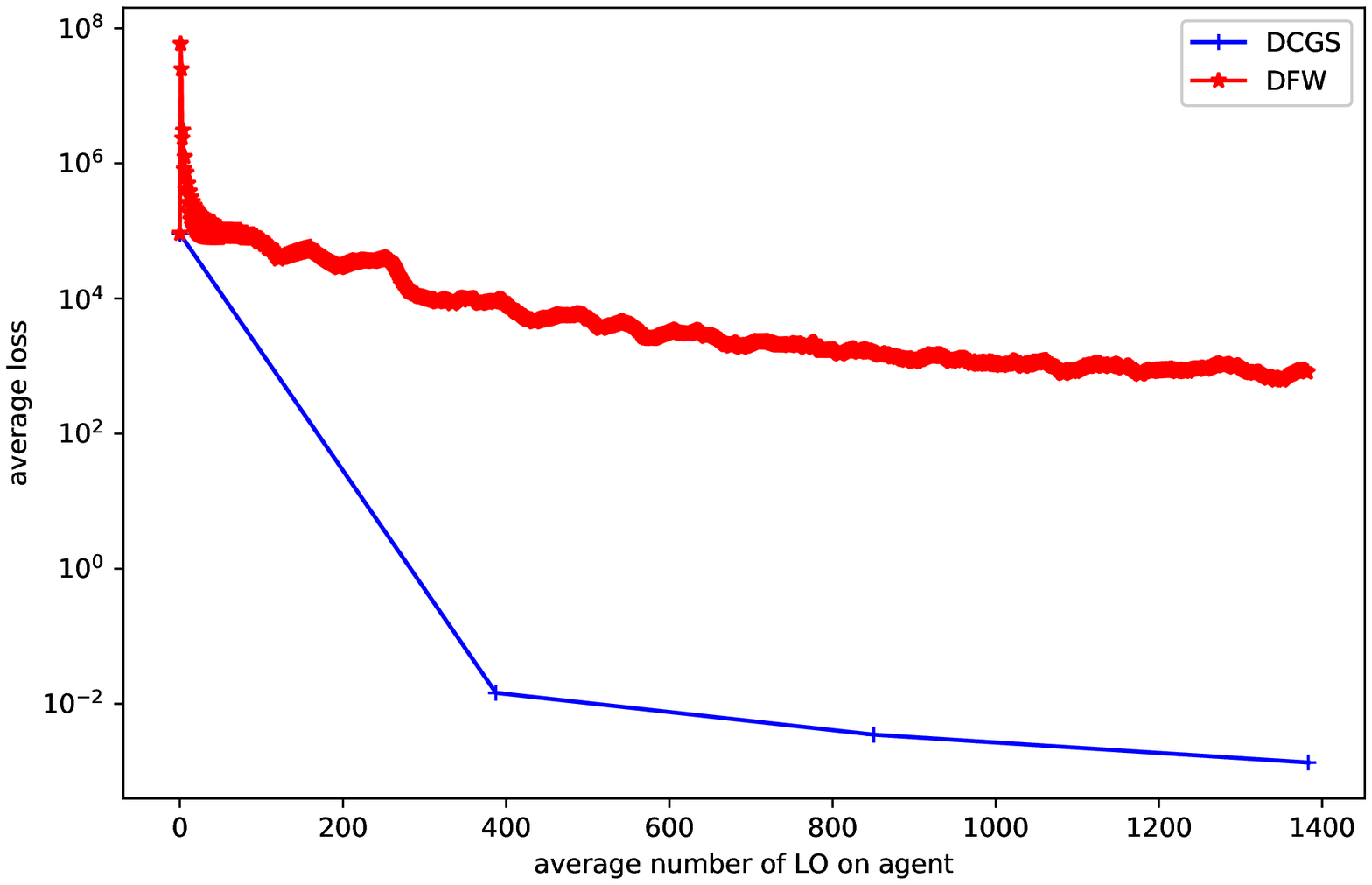}
  \caption{LO vs. Loss}
  \label{fig:synthetic_lo}
\end{subfigure}%
\begin{subfigure}{.24\textwidth}
  \includegraphics[width=\linewidth]{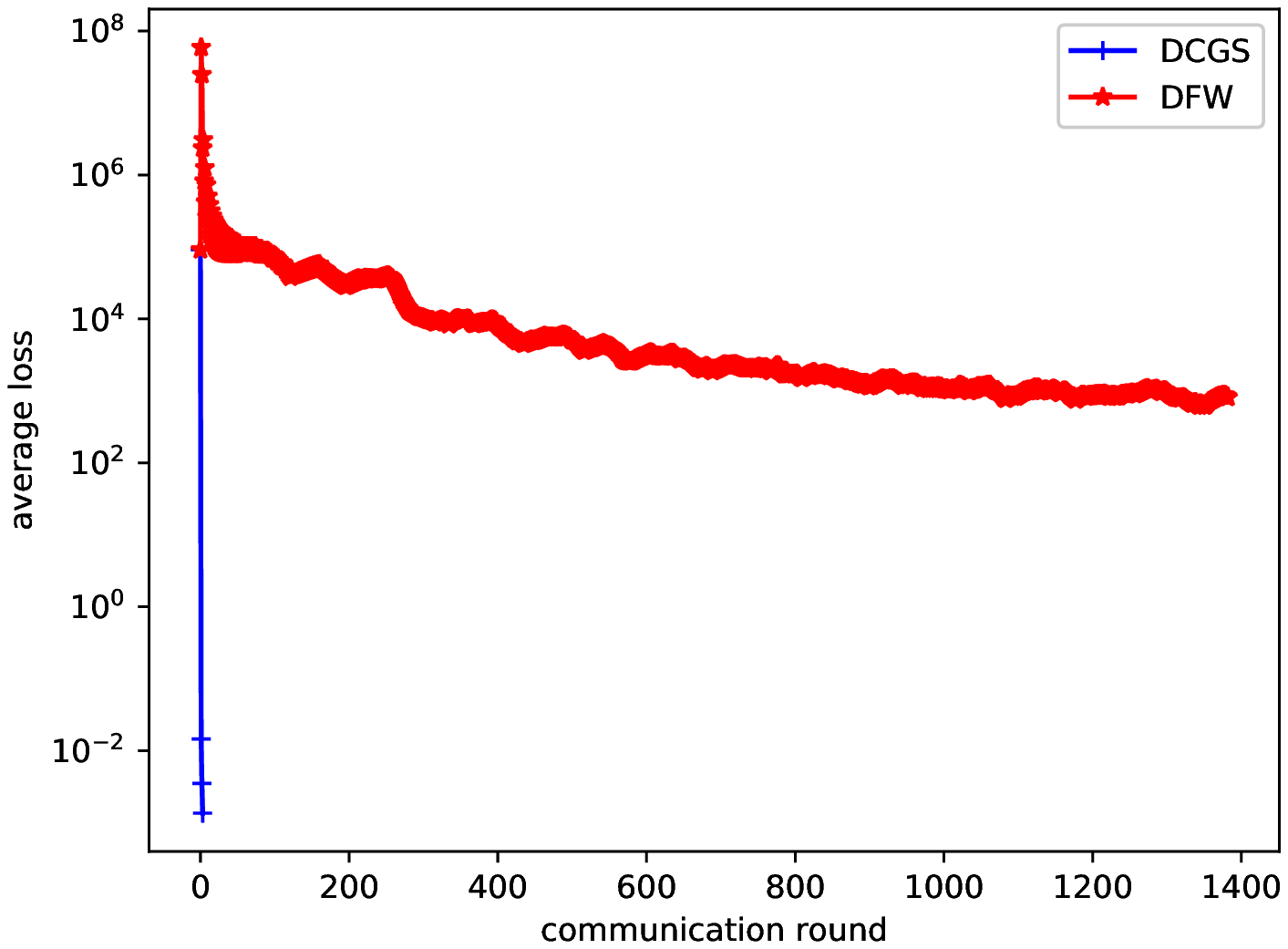}
  \caption{Comm. vs. Loss}
  \label{fig:synthetic_comm}
\end{subfigure}%
\begin{subfigure}{.25\textwidth}
  \includegraphics[width=\linewidth]{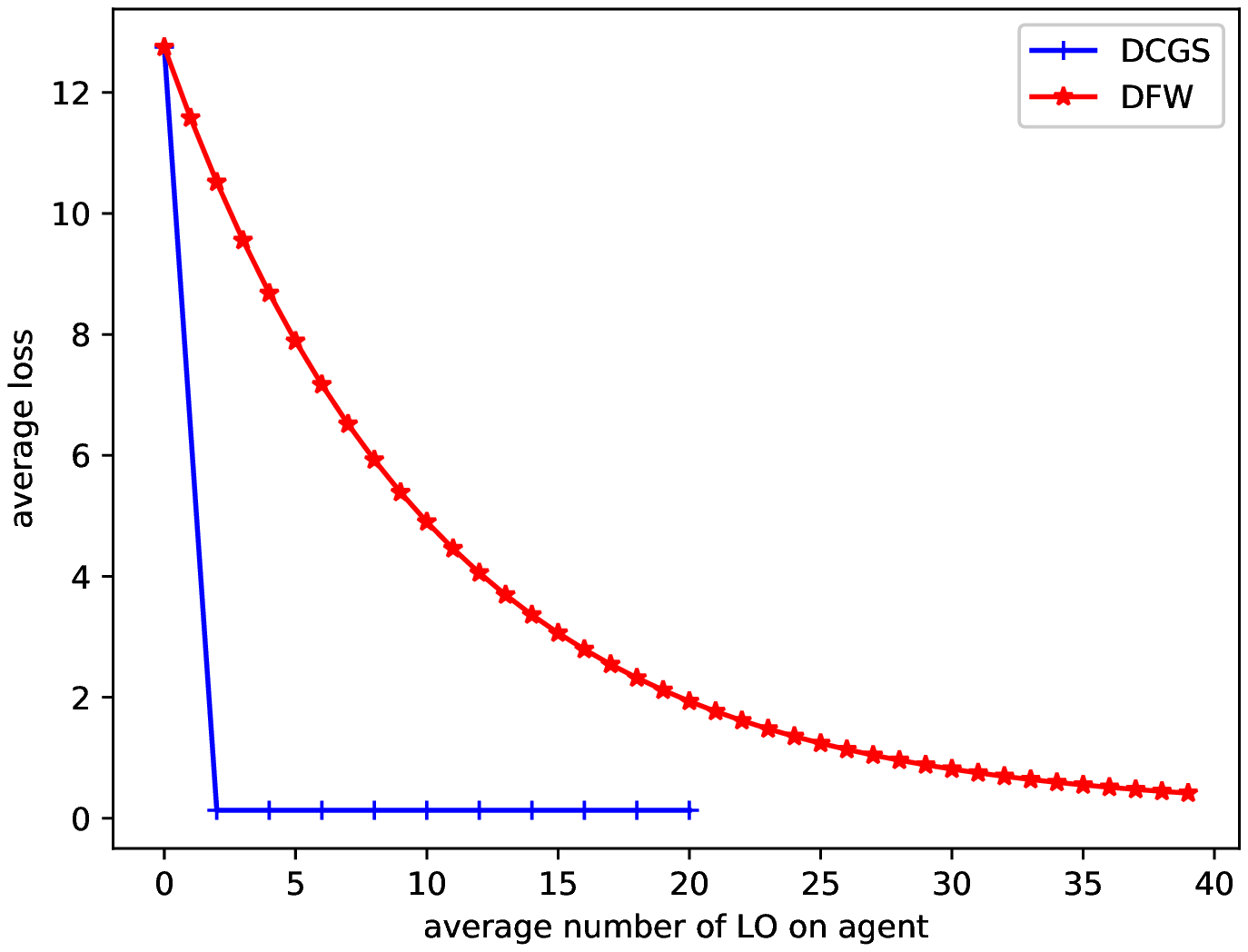}
  \caption{LO vs. Loss}
  \label{fig:real_lo}
\end{subfigure}%
\begin{subfigure}{.25\textwidth}
  \includegraphics[width=\linewidth]{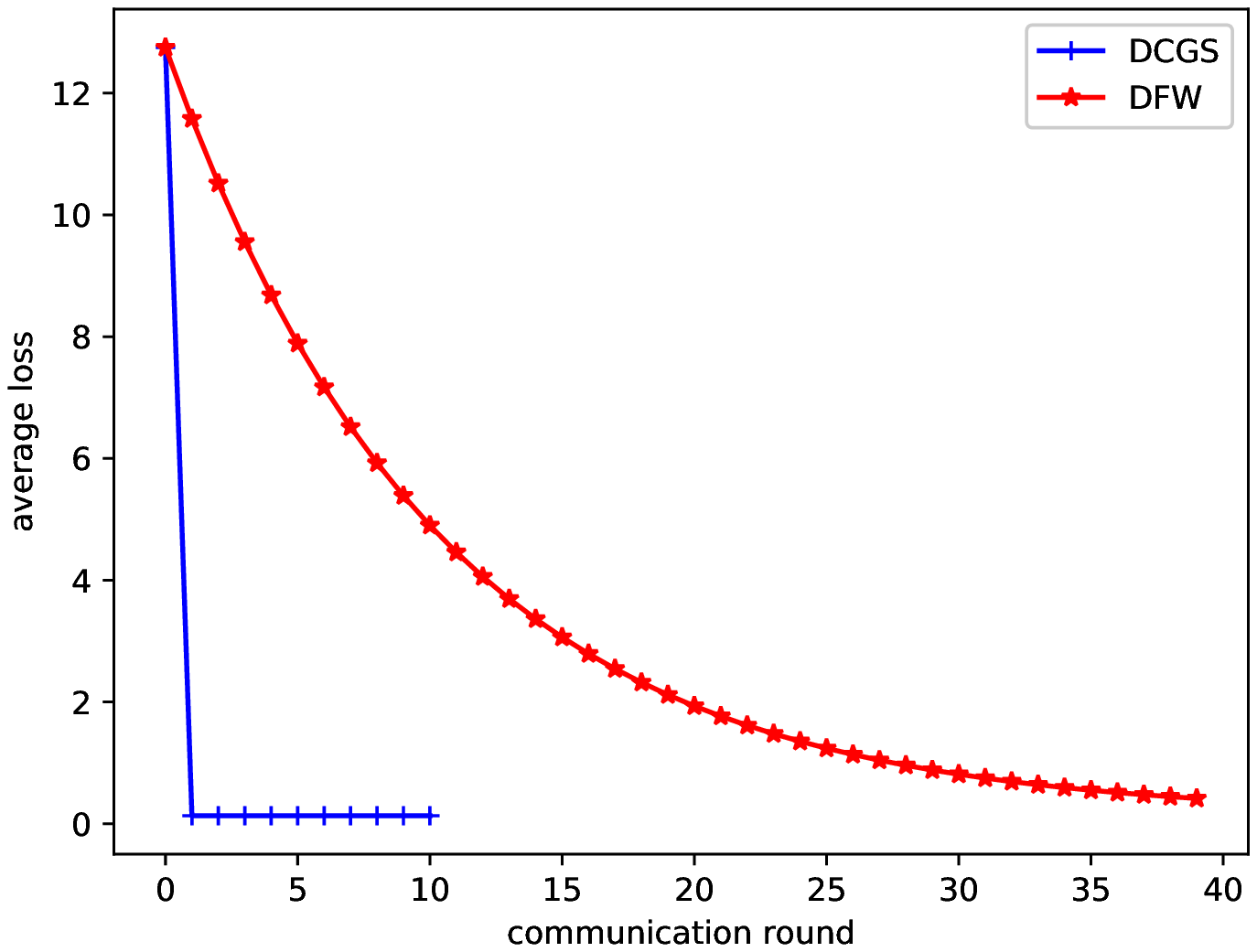}
  \caption{Comm. vs. Loss}
  \label{fig:real_comm}
\end{subfigure}
\caption{Lasso. (a)(b) are for the synthetic data; (c)(d) are for the real data.}
\label{fig:Lasso}
\end{figure}

\textbf{Lasso.}
We compare DCGS and DFW applying to the Lasso problem
 on a synthetic dataset and E2006-tfidf dataset \citep{e2006_dataset}. The Lasso problem is formulated as:
\begin{align}
  \min_{\norm{\theta}_1  \leqslant \rho} \norm{X \theta - Y}^2 & = \sum_{i=1}^n (X_i^{\intercal}\theta - Y_i)^2.
\end{align}
Similar experiment was also conducted in \citep{julien_linear_convergence} which showed linear convergence of FW variant over the polyhedral set. For synthetic data, we generate $n=2000$ samples, with $X_i$ sampled i.i.d.\ from $N(0, I_d)$ and $d=10000$. We
generate $\norm{\theta_{0}} = 100$ with randomly selected $100$ non-zero
entries, and $Y_i = X_i^{\intercal} \theta_{0} + \epsilon_i $ with
$\epsilon_i \sim N(0,1)$. For the real dataset, $d = 150360$ and we randomly draw $n=4000$ samples from the entire E2006-tfidf dataset. We set $\rho = 10^3$ and distribute data evenly across
agents. The results are presented in Figure \ref{fig:Lasso}.
We observed significant performance improvement of DCGS over DFW. For the synthetic data, Figure (\ref{fig:synthetic_lo}, \ref{fig:synthetic_comm}) shows that DCGS with a moderately number of LO, converges to a high-quality solution that has loss by orders of magnitude better than DFW. The gap on communication complexity is even more significant. DFW algorithm takes more than $800$ rounds of communications while DCGS only takes $3$ rounds. We observe similar performance gap on E2006-tfidf dataset in Figure (\ref{fig:real_lo},\ref{fig:real_comm}).

\begin{figure}[h]
\begin{subfigure}{.25\textwidth}
  \includegraphics[width=\linewidth]{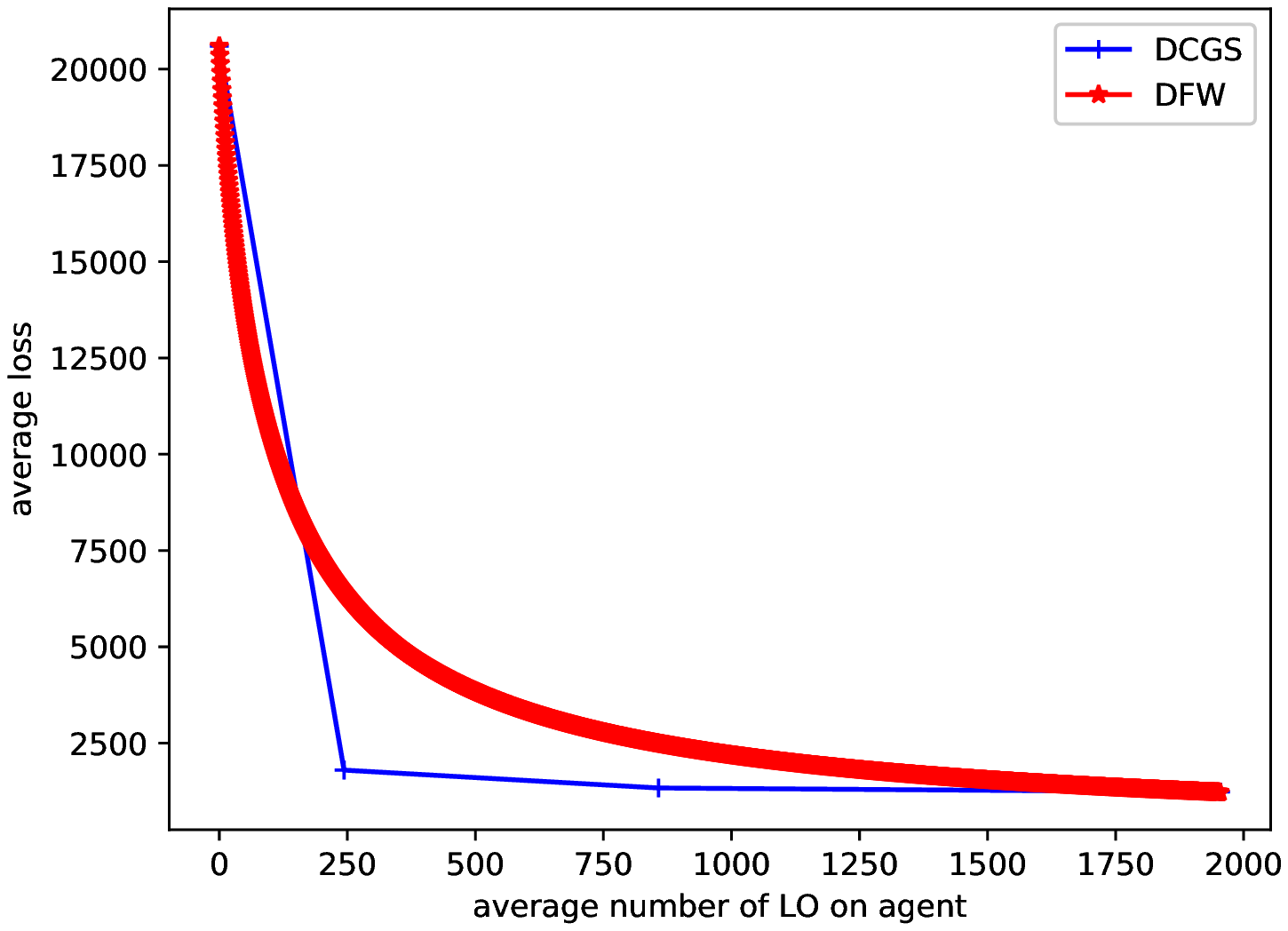}
  \caption{LO vs. Loss}
  \label{fig:synthetic_lo_completion}
\end{subfigure}%
\begin{subfigure}{.25\textwidth}
  \includegraphics[width=\linewidth]{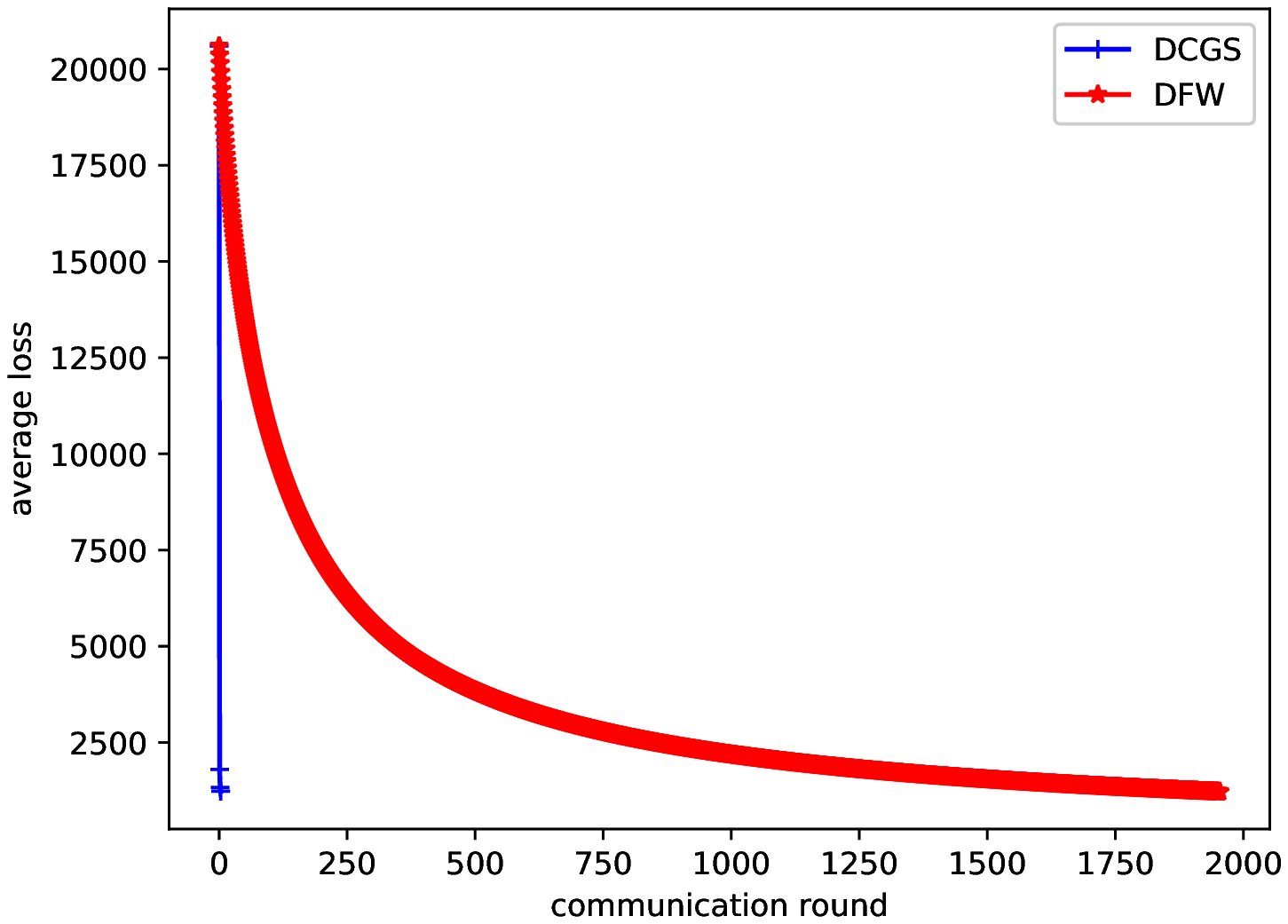}
  \caption{Comm. vs. Loss}
  \label{fig:synthetic_comm_completion}
\end{subfigure}%
\begin{subfigure}{.25\textwidth}
  \includegraphics[width=\linewidth]{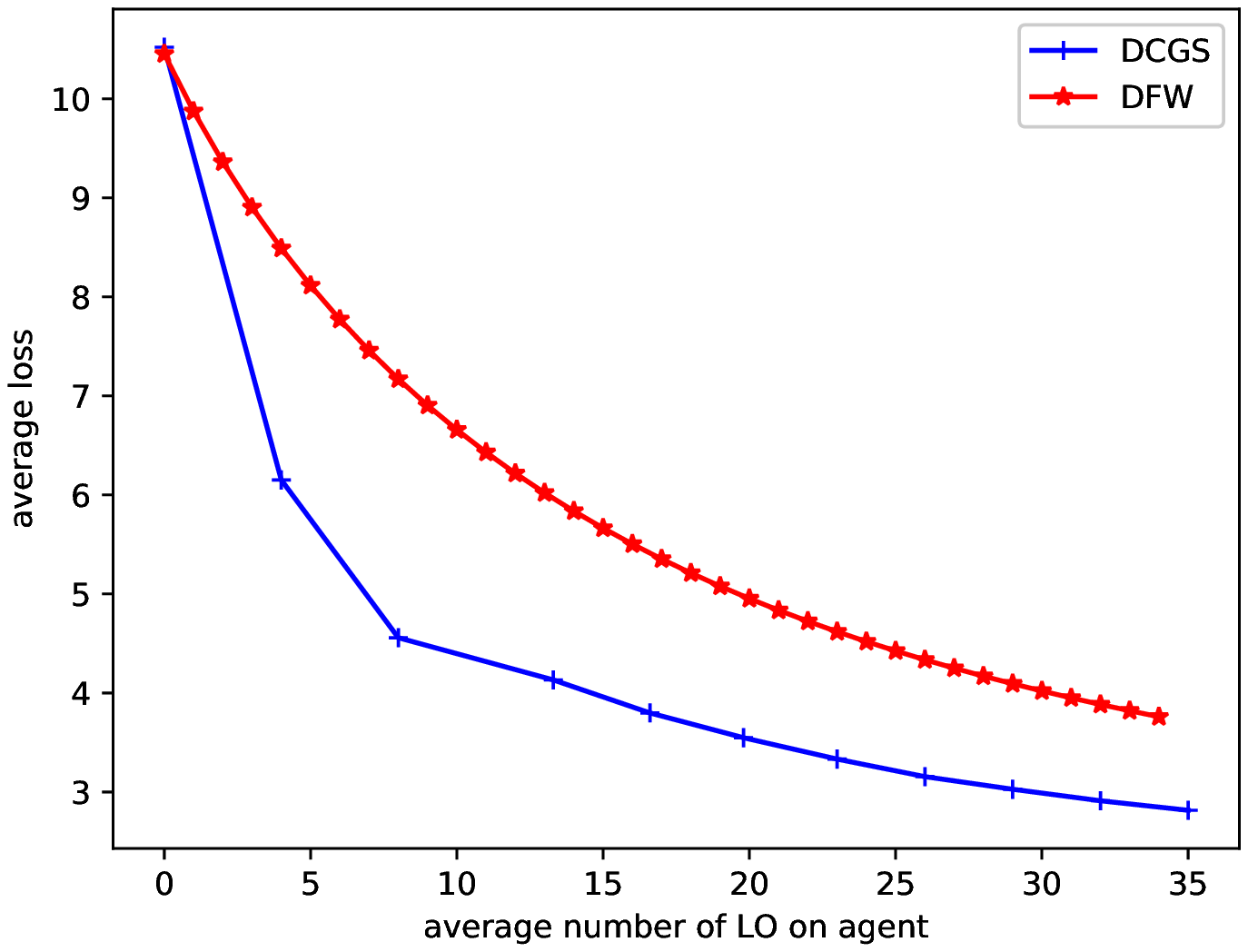}
  \caption{LO vs. Loss}
  \label{fig:real_lo_completion}
\end{subfigure}%
\begin{subfigure}{.25\textwidth}
  \includegraphics[width=\linewidth]{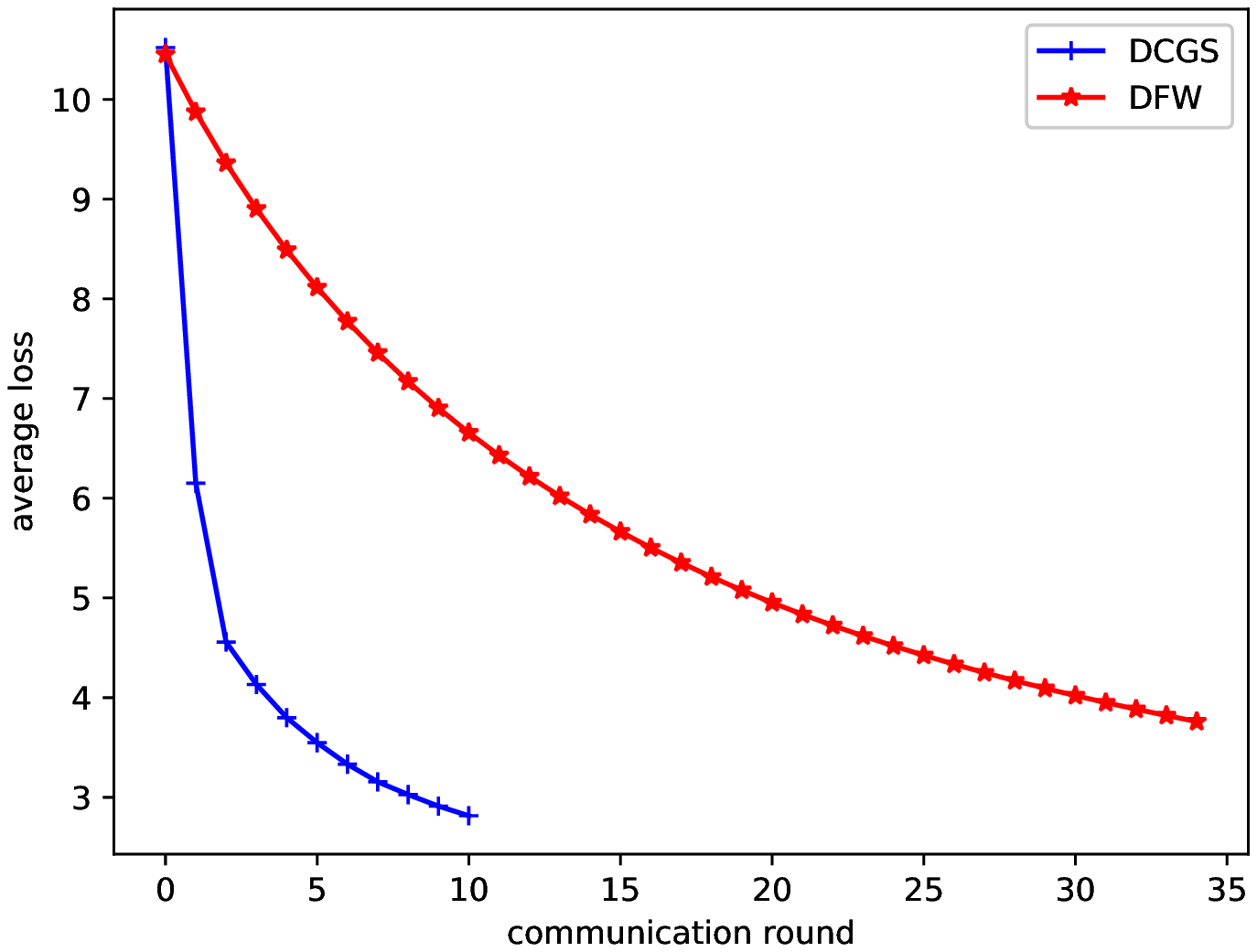}
  \caption{Comm. vs. Loss}
  \label{fig:real_comm_completion}
\end{subfigure}%
\caption{Matrix Completion. (a)(b) are for the synthetic data; (c)(d) are for the real data.}
\label{fig:completion}
\end{figure}

\textbf{Matrix completion.}
We compare DCGS and DFW applying to matrix completion problems on synthetic dataset and MovieLens 100K dataset \citep{movie100k_dataset}. We remark that matrix completion is in fact a communication expensive task. As a toy example, consider a $10^4 \times 10^4$ matrix which takes $800$ MB memory, sending this matrix with 10 MB/s network speed takes $80$ seconds, however an LO on a 4-core machine with Intel(R) Core(TM) i5-6267U CPU @ 2.90GHz processor and 16GB RAM takes less than $2$ seconds. This means for algorithms such as DFW, over $97\%$ of computation time would be waiting for the communication to complete.
We present our simulation results in Figure (\ref{fig:completion}). For the synthetic data set we generate the ground truth matrix $\Theta^* \in \mathbb{R}^{d\times d}$ with $d = 200$ and rank $r =5$. We randomly sample $n=5000$ entries $\{\left(a(i),b(i)\right)\}_{i=1}^n$ and observe
$\Theta^P_{a(i)b(i)} = \Theta^*_{a(i)b(i)} + \epsilon_i $with $\epsilon_i\sim N(0,0.1)$. For MovieLens 100K dataset we want to recover the rating matrix $\Theta^* \in \mathbb{R}^{d_1 \times d_2} $ with $d_1 = 943$, $d_2 = 1682$ and we observed $n=10^5$ ratings.
We set $\rho = 10^4$ and run the algorithms on solving the following objective:
\begin{align}
 \min_{\norm{X}_* \leqslant \rho} \sum_{i=1}^n \left(X_{a(i)b(i)} - \Theta^P_{a(i)b(i)} \right)^2.
\end{align}
For the synthetic data, Figure \ref{fig:synthetic_lo_completion} and \ref{fig:synthetic_comm_completion} show DCGS and DFW need comparable LO to converge to a moderate precision solution, however DFW takes significantly more rounds of communication   ($800 \text{ vs. } 5$). Since communication is the main computation bottle as we discussed above, DCGS would significantly outperform DFW  in terms the actual running time. We observe similar performance gap on MovieLens 100K dataset in Figure (\ref{fig:real_lo_completion},\ref{fig:real_comm_completion}). Our experiment results thus suggest the applicability of DCGS in communication expensive applications.

\section{Conclusions}
In this paper, we propose a communication efficient, distributed projection-free algorithm called DCGS. We show that DCGS is communication efficient under convex and strongly convex setting without restrictive assumptions in existing work, and demonstrate the superiority of DCGS in communication expensive learning tasks such as matrix completion. We also show DCGS can be further improved when the feasible set is polyhedral, which is also validated by our numerical experiments. Future research directions include developing asynchronous DCGS variant and extending DCGS to non-convex settings.

\newpage
\bibliographystyle{elsarticle-harv}
\bibliography{DCGS}

\newpage
\appendix
\rule{\textwidth}{4px}\\
\begin{center}
\textbf{\LARGE Supplemental Material\\
\hfill \\
\Large Communication Efficient Projection Free Algorithm for Distributed Optimization}
\end{center}
\rule{\textwidth}{1px}\\

\textbf{Outline.} In Appendix \ref{convergence} we prove the convergence results for Decentralied Conditional Gradient Sliding. In Appendix \ref{complexity} we establish the communication complexity and the LO complexity of DCGS for the general constraint set. In Appendix \ref{complexity:poly} we establish the communication complexity and the LO complexity of DCGS for the polyhedral constraint set.

\section{Proof of Main Theorem}\label{convergence}
\subsection{Proof of Theorem \ref{thrm:smooth}}
In this subsection we prove the convergence result for smooth and convex objective. Denote $Y = \mathbb{R}^{md}$
 and recall the saddle point problem defined in (\ref{bilinear_saddle_formulation}
),for $\mathbf{z} = (\mathbf{x}, \mathbf{y}) , \overline{\mathbf{z}} = (\overline{\mathbf{x}}, \overline{\mathbf{y}}) \in X^m \times Y$,
we define the primal-dual gap function to be:
\begin{align}\label{vanilla_gap_function}
  Q(\z, \overline{\z}) = F(\mathbf{x}) + \inner{\mathbf{Lx}}{\mathbf{\overline{y}}} - F(\overline{\mathbf{x}}) - \inner{\mathbf{L\overline{x}}}{\mathbf{y}}
\end{align}
Note that if $\z^* = (\x^*, \y^*)$ is a saddle point to (\ref{bilinear_saddle_formulation}
), then $Q(\z^*, \z) \leqslant 0 $ and $Q(\z, \z^*)\geqslant 0 $ for any $\z \in X^m \times Y$. It is then natural to measure the quality of a solution $\z$ to problem (\ref{bilinear_saddle_formulation}
) by $\sup_{\overline{\z} \in X^m \times Y} Q(\z^*, \overline{\z})$.
To handle unboundedness of $Y$ here, we define the modified gap function to be:
\begin{align}\label{gap_function}
  g_{\mathbf{Y}} (\mathbf{s},\mathbf{z}) = \sup_{\overline{\mathbf{y}} \in Y} Q(z; (\mathbf{x}^*, \overline{\mathbf{y}})) - \inner{\mathbf{s}}{\overline{\mathbf{y}}}
\end{align}
 In fact, we have the following proposition.
\begin{prop}[\citep{lan_comm_efficient}]
 If we have $g_{Y} (\mathbf{s},\mathbf{z}) \leqslant \epsilon $ for any $\epsilon >0$, then we must have $ \mathbf{Lx} = \mathbf{s}$ and $F(\mathbf{x}) - F(\mathbf{x}^*) \leqslant \epsilon$ .
\end{prop}
  this claim is straightforward to eastablish by following the definition of $g_{Y} (\mathbf{s},\mathbf{z}) $ and a proof by contradiction argument.

By construction of $y_i^k$ in Algorithm \ref{alg:dcgs}
 we know that:
\begin{align}\label{proj_lemma_y}
 \inner{v_i^k}{y_i^k - y_i} = \frac{\tau_k}{2} \left[ \norm{y_i - y_i^{k-1}}^2 - \norm{y_i - y_i^k}^2 - \norm{y_i^{k-1} - y_i^k}^2  \right]
\end{align}
Similar inequality could be established for $x_i^k$. Recall in Line \ref{lo_subproblem} of Algorithm \ref{alg:dcgs} we are solving the subproblem
\begin{align}
\min_{x_i \in X }\phi_i^k (x_i) = \inner{w_i^k}{x_i} + f_i(x_i) + \frac{\eta_k}{2} \norm{x_i - x_i^{k-1}}^2
\end{align}
with the ouput $x_i^k$ satisfying inequality $\inner{\phi_i^k(x_i^k)}{x_i^k - x_i^{k-1}} \leqslant e_i^k$.  Since $\phi_i^k(x_i)$ is strongly convex, we have: $\phi_i^k(x_i) - \phi_i^k(x_i^k) \geqslant \inner{\nabla \phi_i^k(x_i^k)}{x_i - x_i^k} +  \frac{\eta_k}{2} \norm{x_i -x_i^k}$. Combine this two inequalities with some algebraic rearrangements yields the following:
\begin{align}\label{proj_lemma_x}
 \inner{w_i^k}{x_i^k - x_i} + f_i(x_i^k) - f_i(x_i) \leqslant
 \frac{\eta_k}{2}\left[ \norm{ x_i - x_i^{k-1}}^2 - \norm{x_i - x_i^k}^2 -\norm{x_i^{k-1} - x_i^k}^2\right] + e_i^k
\end{align}
Summing up the previous two inequalities and using the definition of $Q(\cdot, \cdot), w_i^k, v_i^k$ we have:
\begin{align}
 Q(\mathbf{z^k}, \mathbf{z}) & = F(\mathbf{x^k}) - F(\mathbf{x}) +
 \inner{\mathbf{Lx^k}}{\mathbf{y}} - \inner{\L \x}{\y^k} ]\nonumber \\
  & = \inner{\L (\x^k - \tilde{\x}^k)}{\y - \y^k} + \frac{\eta_k}{2} \left[ \norm{\x^{k-1}- \x}^2 - \norm{ \x^k -\x}^2 - \norm{\x^{k-1} - \x^k}^2\right] \nonumber  \\
  & \quad + \frac{\tau_k}{2} \left[\norm{\y - \y^{k-1}}^2 - \norm{\y - \y^k}^2 - \norm{\y^k - \y^{k-1}}^2 \right] + \sum_{i=1}^m e_i^k \label{initial_bound}
\end{align}
We define the right hand side of previous equation by $\Delta_k$, and we are going to handle the weighted sum of the first three terms in $\Delta_k$ seperately. For the first term:
\begin{align}
\sum_{k=1}^N \theta_k \inner{\L (\x^k - \tilde{\x}^k)}{\y - \y^k} & =
\sum_{k=1}^N \theta_k \left[\inner{\L (\x^k - \x^{k-1})}{\y - \y^k}
- \alpha_k  \inner{\L (\x^{k-1} - \x^{k-2})}{\y - \y^k}\right] \label{extrapolation}\\
& = \sum_{k=1}^N \left[ \theta_k \inner{\L (\x^k - \x^{k-1})}{\y - \y^k}
- \alpha_k \theta_k  \inner{\L (\x^{k-1} - \x^{k-2})}{\y - \y^{k-1}} \right] \nonumber \\
& \quad + \sum_{k=1}^N \alpha_k \theta_k  \inner{\L (\x^{k-1} - \x^{k-2})}{\y^{k-1} - \y^{k}}  \nonumber \\
& = \sum_{k=1}^N \left[ \theta_k \inner{\L (\x^k - \x^{k-1})}{\y - \y^k}
- \theta_{k-1}  \inner{\L (\x^{k-1} - \x^{k-2})}{\y - \y^{k-1}} \right]  \nonumber \\
& \quad + \sum_{k=1}^N \alpha_k \theta_k  \inner{\L (\x^{k-1} - \x^{k-2})}{\y^{k-1} - \y^{k}}   \label{par_choice} \\
& = \theta_N \inner{\L (\x^N - \x^{N-1})}{\y - \y^N} \nonumber \\
& \quad + \sum_{k=1}^N \alpha_k \theta_k  \inner{\L (\x^{k-1} - \x^{k-2})}{\y^{k-1} - \y^{k}}
\label{tele3}
\end{align}
where (\ref{extrapolation}) follows from definition of $\tilde{x_i}^k$. In (\ref{par_choice}) we use condition $\alpha_k \theta_k = \theta_{k-1}$ which follows from our parameters setting. (\ref{tele3}) comes from telescoping the first summation in (\ref{par_choice}) and the condition that $\x^0 = \x^{-1}$.
We could bound the weighted sum of the second term in (\ref{initial_bound}) by:
\begin{align}
\sum_{k=1}^N \frac{\theta_k \eta_k}{2} [  \norm{\x^{k-1}- \x}^2 - \norm{ \x^k -\x}^2 & -  \norm{\x^{k-1} - \x^k}^2 ]  =
\frac{\theta_1 \eta_1}{2} \norm{\x^0 - \x}^2  - \frac{\theta_N \eta_N}{2} \norm{\x^N - \x}^2 \nonumber \\
& + \sum_{k=1}^{N-1}
\frac{\theta_{k+1}\eta_{k+1} - \theta_k \eta_k}{2} \norm{\x^k - \x}^2
- \sum_{k=1}^N \frac{\theta_k \eta_k}{2} \norm{\x^{k-1} - \x^k}^2 \nonumber \\
& \leqslant \frac{\theta_1 \eta_1}{2} \norm{\x^0 - \x}^2  - \frac{\theta_N \eta_N}{2} \norm{\x^N - \x}^2 \nonumber \\
& - \sum_{k=1}^N \frac{\theta_k \eta_k}{2} \norm{\x^{k-1} - \x^k}^2 \label{tele2}
\end{align}
where in (\ref{tele2}) we use the condition $\theta_{k+1}\eta_{k+1} \leqslant \theta_k \eta_k $ which follows from our paramters setting. Similarly  we can bound the weighted sum of the third term in (\ref{initial_bound}):
\begin{align}
\sum_{k=1}^N \frac{\theta_k \eta_k}{2} [  \norm{\y^{k-1}- \x}^2 - \norm{ \y^k -\y}^2 & -  \norm{\y^{k-1} - \y^k}^2 ]  =
\frac{\theta_1 \eta_1}{2} \norm{\y^0 - \y}^2  - \frac{\theta_N \eta_N}{2} \norm{\y^N - \y}^2 \nonumber \\
& + \sum_{k=1}^{N-1}
\frac{\theta_{k+1}\eta_{k+1} - \theta_k \eta_k}{2} \norm{\y^k - \y}^2
- \sum_{k=1}^N \frac{\theta_k \eta_k}{2} \norm{\y^{k-1} - \y^k}^2 \nonumber \\
& \leqslant \frac{\theta_1 \eta_1}{2} \norm{\y^0 - \y}^2  - \frac{\theta_N \eta_N}{2} \norm{\y^N - \y}^2 \nonumber \\
& - \sum_{k=1}^N \frac{\theta_k \eta_k}{2} \norm{\y^{k-1} - \y^k}^2 \label{tele1}
\end{align}
Now sum up (\ref{tele3}),(\ref{tele2}),(\ref{tele1}) we get:
\begin{align}
\sum_{k=1}^N \theta_k \Delta_k &
   \leqslant \theta_N \inner{\L (\x^N - \x^{N-1})}{\y - \y^N} \nonumber \\
  & \quad + \sum_{k=1}^N \theta_k [ \alpha_k \inner{\L (\x^{k-1} - \x^{k-2})}{\y^{k-1} - \y^{k}} - \frac{\eta_k}{2} \norm{\x^{k-1} - \x^k}^2 - \frac{\tau_k}{2} \norm{\y^{k-1}-\y^k}^2] \nonumber \\
  & \quad + \frac{\theta_1 \eta_1}{2} \norm{\x^0 - \x}^2
   - \frac{\theta_N \eta_N}{2} \norm{\x^N - \x}^2 \nonumber \\
  & \quad + \frac{\theta_1 \tau_1}{2} \norm{\y^1 - \y}^2
   - \frac{\theta_N \tau_N}{2} \norm{\y^N - \y}^2
   + \sum_{k=1}^N \sum_{i=1}^m \theta_k e_i^k \label{initial_cond}
\end{align}
We further rewrite the first summation term in (\ref{initial_cond}) as the following:
\begin{align}
 & \sum_{k=1}^N \theta_k [ \alpha_k \inner{\L (\x^{k-1} - \x^{k-2})}{\y^{k-1} - \y^{k}} - \frac{\eta_k}{2} \norm{\x^{k-1} - \x^k}^2 - \frac{\tau_k}{2} \norm{\y^{k-1}-\y^k}^2] \nonumber \\
 = &  \sum_{k=2}^N \left[\theta_k \alpha_k \inner{\L (\x^{k-1} - \x^{k-2})}{\y^{k-1} - \y^{k}}
-  \frac{\theta_{k-1}\eta_{k-1}}{2} \norm{\x^{k-2} - \x^{k-1}}^2 -  \frac{\theta_k \tau_k}{2} \norm{\y^{k-1}-\y^k}^2 \right ]  \nonumber \\ & \quad - \frac{\theta_N \eta_N}{2} \norm{\x^{N-1} - \x^N}^2 - \frac{\eta_1 \tau_1}{2} \norm{\y^0 - \y^1}^2 \label{rewrite}
\end{align}
Combine (\ref{initial_cond}) and (\ref{rewrite}) we obtain:
\begin{align}
\sum_{k=1}^N \theta_k \Delta_k & \leqslant \theta_N \inner{\L (\x^N - \x^{N-1})}{\y - \y^N} - \frac{\theta_N \eta_N}{2} \norm{\x^{N-1} - \x^N}^2 - \frac{\theta_1 \tau_1}{2} \norm{y^0 - y^1}^2 \nonumber \\
&  \quad+ \sum_{k=2}^N \left[\theta_k \alpha_k \inner{\L (\x^{k-1} - \x^{k-2})}{\y^{k-1} - \y^{k}}
-  \frac{\theta_{k-1}\eta_{k-1}}{2} \norm{\x^{k-2} - \x^{k-1}}^2 -  \frac{\theta_k \tau_k}{2} \norm{\y^{k-1}-\y^k}^2 \right ] \nonumber \\
& \quad + \frac{\theta_1 \eta_1}{2} \norm{\x^0 - \x}^2 - \frac{\theta_N \eta_N}{2} \norm{\x^N - \x}^2 \nonumber \\
& \quad + \frac{\theta_1 \tau_1}{2} \norm{\y^0 -\y}^2 - \frac{\theta_N \tau_N}{2}\norm{\y^N - \y}^2 + \sum_{k=1}^N \sum_{i=1}^m \theta_k e_i^k
\end{align}
The summation in the second line could be in fact upper bounded by $0$ as the following:
\begin{align}
  & \sum_{k=2}^N \left[\theta_k \alpha_k \inner{\L (\x^{k-1} - \x^{k-2})}{\y^{k-1} - \y^{k}}
  -  \frac{\theta_{k-1}\eta_{k-1}}{2} \norm{\x^{k-2} - \x^{k-1}}^2 -  \frac{\theta_k \tau_k}{2} \norm{\y^{k-1}-\y^k}^2 \right ] \nonumber \\
  \leqslant & \sum_{k=2}^N \left[\theta_k \alpha_k \norm{\L} \norm{\x^{k-1} - \x^{k-2}}\norm{\y^{k-1} - \y^k }
  -  \frac{\theta_{k-1}\eta_{k-1}}{2} \norm{\x^{k-2} - \x^{k-1}}^2 -  \frac{\theta_k \tau_k}{2} \norm{\y^{k-1}-\y^k}^2 \right ] \nonumber\\
  \leqslant & \sum_{k=2}^N \left(\frac{\norm{\L}^2}{2\tau_k \theta_k}
  - \frac{\theta_{k-1}\eta_{k-1}}{2} \right) \norm{\x^{k-2} - \x^{k-1}}^2 \label{young} \\
  = & \sum_{k=2}^N \left(\frac{\theta_{k-1}\alpha_k \norm{\L}^2}{2\tau_k}
  - \frac{\theta_{k-1}\eta_{k-1}}{2} \right) \norm{\x^{k-2} - \x^{k-1}}^2 \leqslant 0\label{par_choice2}
\end{align}
where in (\ref{young}) we use the Young's inequality, in \ref{par_choice2} we use the condition $\theta_k \alpha_k = \theta_{k-1}, \alpha_k \norm{\L}^2 \leqslant \tau_k \eta_{k-1}$ which follows from our parameter setting.
In summary we get:
\begin{align}
  \sum_{k=1}^N \theta_k Q(\z^k, \z) = \sum_{k=1}^N \theta_k \Delta_k & \leqslant \theta_N \inner{\L (\x^N - \x^{N-1})}{\y - \y^N} - \frac{\theta_N \eta_N}{2} \norm{\x^{N-1} - \x^N}^2 \nonumber \\
  & \quad + \frac{\theta_1 \eta_1}{2} \norm{\x^0 - \x}^2 - \frac{\theta_N \eta_N}{2} \norm{\x^N - \x}^2 \nonumber \\
  & \quad + \frac{\theta_1 \tau_1}{2} \norm{\y^0 -\y}^2 - \frac{\theta_N \tau_N}{2}\norm{\y^N - \y}^2 +
  \sum_{k=1}^N \sum_{i=1}^m \theta_k e_i^k\label{crude}
\end{align}
Our next objective is to bound the right hand side of (\ref{crude}) as a linear function on $\y$. Collecting all the linear term of $\y$ after some rearrangement, we get:
\begin{align}
  \sum_{k=1}^N \theta_k Q(\z^k, \z)  & \leqslant
  \theta_N \inner{\y^N}{\L(\x^{N-1} - \x^N)} - \frac{\theta_N \eta_N}{2}
  \norm{ \x^{N-1} - \x^N} - \frac{\theta_N \tau_N}{2} \norm{\y^N}^2 \nonumber \\
  & \quad + \frac{\theta_1 \eta_1}{2} \norm{\x^0 - \x}^2 + \frac{\theta_1 \tau_1}{2} \norm{\y^0}^2 \nonumber \\
  & \quad + \inner{\y}{\theta_N \L(\x^N - \x^{N-1}) + \theta_1 \tau_1 (\y^N - \y^0)} +\sum_{k=1}^N \sum_{i=1}^m \theta_k e_i^k
  \nonumber \\
  & \leqslant \theta_N \norm{\L} \norm{\x^{N-1} - \x^N} \norm{\y^N} - \frac{\theta_N \eta_N}{2}
  \norm{ \x^{N-1} - \x^N} - \frac{\theta_N \tau_N}{2} \norm{\y^N}^2 \nonumber \\
  & \quad + \frac{\theta_1 \eta_1}{2} \norm{\x^0 - \x}^2 + \frac{\theta_1 \tau_1}{2} \norm{\y^0}^2 \nonumber \\
  & \quad+ \inner{\y}{\theta_N \L(\x^N - \x^{N-1}) + \theta_1 \tau_1 (\y^N - \y^0)} \nonumber + \sum_{k=1}^N \sum_{i=1}^m \theta_k e_i^k \nonumber \\
  & \leqslant \left(\frac{\theta_N \norm{\L}^2}{2\eta_N} - \frac{\theta_1 \tau_1}{2}  \right) \norm{\y^N}^2
  + \frac{\theta_1 \eta_1}{2} \norm{\x^0 - \x}^2 + \frac{\theta_1 \tau_1}{2} \norm{\y^0}^2  \nonumber \\
  &\quad + \inner{\y}{\theta_N \L(\x^N - \x^{N-1}) + \theta_1 \tau_1 (\y^N - \y^0)} +\sum_{k=1}^N \sum_{i=1}^m \theta_k e_i^k \label{young1}\\
  & \leqslant \frac{\theta_1 \eta_1}{2} \norm{\x^0 - \x}^2 + \frac{\theta_1 \tau_1}{2} \norm{\y^0}^2  \nonumber \\
  & \quad + \inner{\y}{\theta_N \L(\x^N - \x^{N-1}) + \theta_1 \tau_1 (\y^N - \y^0)} + \sum_{k=1}^N \sum_{i=1}^m \theta_k e_i^k\label{par_choice3}
\end{align}
where in (\ref{young1}) we use the Young's inequality, in (\ref{par_choice3}) we used the condition that $\theta_N \norm{\L}^2 \leqslant \theta_1 \tau_1 \eta_N$ which is satisfied by our parameters. Let us define $\s^N = \theta_N \L(\x^N - \x^{N-1}) + \theta_1 \tau_1 (\y^N - \y^0) $, then we have shown that:
\begin{align}
\sum_{k=1}^N \theta_k Q(\z^k, \z)  & \leqslant
\frac{\theta_1 \eta_1}{2} \norm{\x^0 - \x}^2 + \frac{\theta_1 \tau_1}{2} \norm{\y^0}^2 +\sum_{k=1}^N \sum_{i=1}^m \theta_k e_i^k + \inner{\y}{\s^N}
\end{align}
Choosing $\z =  (\x^*, \y)$ in the left hand side, and using the convexity of $Q(\cdot, \z)$, we immediately have:
\begin{align}
 Q(\overline{\z}^N; (\x^*, \y)) \leqslant (\sum_{k=1}^N \theta_k)^{-1}
 \left( \frac{\theta_1 \eta_1}{2} \norm{\x^0 - \x^*}^2 + \frac{\theta_1 \tau_1}{2} \norm{\y^0}^2 +\sum_{k=1}^N \sum_{i=1}^m \theta_k e_i^k + \inner{\y}{\s^N} \right)
\end{align}
From the definition of $g_{\mathbf{Y}}(\s,\z)$, and define $\overline{\s}^N = (\sum_{k=1}^N \theta_k)^{-1} \s^N$, we have:
\begin{align}
  g_{\mathbf{Y}}(\overline{\s}^N, \overline{\z}^N) \leqslant (\sum_{k=1}^N \theta_k)^{-1}
  \left( \frac{\theta_1 \eta_1}{2} \norm{\x^0 - \x^*}^2 + \frac{\theta_1 \tau_1}{2} \norm{\y^0}^2  +\sum_{k=1}^N \sum_{i=1}^m \theta_k e_i^k \right)
\end{align}
which then implies $ \L \overline{\x}^N  = \overline{\s}^N $ and:
\begin{align}
  F(\overline{\x}^N) - F(\x^*) & \leqslant (\sum_{k=1}^N \theta_k)^{-1}
\left( \frac{\theta_1 \eta_1}{2} \norm{\x^0 - \x^*}^2 + \frac{\theta_1 \tau_1}{2} \norm{\y^0}^2  +\sum_{k=1}^N \sum_{i=1}^m \theta_k e_i^k\right ) \nonumber
\end{align}
Now plug in choice of $\alpha_1 = \theta_k =1 , \eta_k = 2\norm{\L}, \tau_k = \norm{\L} $ and $e_i^k = \frac{\norm{\L} \max(\norm{\x^0 - \x^*}, \norm{\y^0}^2)}{mN}$ yields our convergence result.

\subsection{Proof of Theorem \ref{thrm:sc}}
In this subsection we prove the convergence result for the smooth and strongly convex objective. For strongly convex $f_i$, again by the update for $x_i^k$ in Algorithm \ref{alg:dcgs}, we have $\inner{ \nabla \phi_i^k (x_i^k)}{x_i^k - x_i^{k-1}} \leqslant e_i^k$. Since $\phi_i^k(x_i)$ is strongly convex we have: $\phi_i^k(x_i) - \phi_i^k(x_i^k) \geqslant \inner{\nabla \phi_i^k(x_i^k)}{x_i - x_i^k} + (u + \frac{\eta_k}{2}) \norm{x_i -x_i^k}$. Combine this two inequalities with some algebraic rearrangements we get:
\begin{align}
  \inner{w_i^k}{x_i^k - x_i} + f_i(x_i^k) - f_i(x_i) \leqslant
  \frac{\eta_k}{2} \norm{ x_i - x_i^{k-1}}^2 - (\frac{\eta_k}{2} + u) \norm{x_i - x_i^k}^2 - \frac{\eta_k}{2}\norm{x_i^{k-1} - x_i^k}^2 + e_i^k \nonumber
\end{align}
Note we still have have (\ref{proj_lemma_y}) since the update for $y_i^k$ does not change. Following the same argument of (\ref{initial_bound}) we have:
\begin{align}
  Q(\mathbf{z^k}, \mathbf{z}) & = F(\mathbf{x^k}) - F(\mathbf{x}) +
  \inner{\mathbf{Lx^k}}{\mathbf{y}} - \inner{\L \x}{\y^k} \nonumber \\
   & = \inner{\L (\x^k - \tilde{\x}^k)}{\y - \y^k} + \left(\frac{\eta_k}{2} \norm{\x^{k-1}- \x}^2 -  (\frac{\eta_k}{2} + u )\norm{ \x^k -\x}^2 - \frac{\eta_k}{2}\norm{\x^{k-1} - \x^k}^2 \right) \nonumber  \\
   & \quad + \frac{\tau_k}{2} \left[\norm{\y - \y^{k-1}}^2 - \norm{\y - \y^k}^2 - \norm{\y^k - \y^{k-1}}^2 \right] + \sum_{i=1}^m e_i^k \label{initial_bound_2}
\end{align}
still we can bound the weighted sum the of first term as in (\ref{tele3}) and the weighted sum of third term as in (\ref{tele1}), the paramter condition required by establishing them is still satisfied. Handling the second term is also essentially the same:
\begin{align}
& \sum_{k=1}^N  \theta_k \left(\frac{\eta_k}{2} \norm{\x^{k-1}- \x}^2 -  (\frac{\eta_k}{2} + u )\norm{ \x^k -\x}^2 - \frac{\eta_k}{2}\norm{\x^{k-1} - \x^k}^2 \right)   \nonumber \\
= & \frac{\theta_1 \eta_1}{2} \norm{\x^0 - \x}^2 + \sum_{k=1}^{N-1} \frac{\theta_{k+1} \eta_{k+1} - \theta_k(\eta_k + u)}{2} \norm{\x^k -\x}^2 - \frac{\theta_N (\eta_N + u)}{2} \norm{\x^N - \x}^2 \nonumber \\
& - \sum_{k=1}^N \frac{\theta_k \eta_k}{2} \norm{\x^{k-1} - \x^k}^2 \nonumber \\
\leqslant & \frac{\theta_1 \eta_1}{2} \norm{\x^0 - \x}^2 - \frac{\theta_N (\eta_N + u)}{2} \norm{\x^N - \x}^2 - \sum_{k=1}^N \frac{\theta_k \eta_k}{2} \norm{\x^{k-1} - \x^k}^2 \label{new_tele2}
\end{align}
where (\ref{new_tele2}) comes from the condition $\theta_{k+1} \eta_{k+1} \leqslant \theta_k(\eta_k + u) $ which follows from our parameter setting.
Now add up (\ref{tele3}), (\ref{tele1}), (\ref{new_tele2}) and combine with (\ref{initial_bound_2}) we have the following:
\begin{align}
  \sum_{k=1}^N \theta_k Q(\z^k, \z) = \sum_{k=1}^N \theta_k \Delta_k & \leqslant \theta_N \inner{\L (\x^N - \x^{N-1})}{\y - \y^N} - \frac{\theta_N \eta_N}{2} \norm{\x^{N-1} - \x^N}^2 \nonumber \\
  & \quad + \frac{\theta_1 \eta_1}{2} \norm{\x^0 - \x}^2 - \frac{\theta_N (\eta_N + u)}{2} \norm{\x^N - x}^2  \nonumber \\
  & \quad + \frac{\theta_1 \tau_1}{2} \norm{\y^0 -\y}^2 - \frac{\theta_N \tau_N}{2}\norm{\y^N - \y}^2 +
  \sum_{k=1}^N \sum_{i=1}^m \theta_k e_i^k \label{final_bound2}
\end{align}
We now need to bound the right hand side of (\ref{final_bound2}) by a linear function on $\y$. This is exactly the same as in
 establishing (\ref{par_choice3}) and hence we omit the tedious detail. In summary, we get the following bound that is identical to (\ref{par_choice3}):
 \begin{align}
  \sum_{k=1}^N \theta_k Q(\z^k, \z) & \leqslant \frac{\theta_1 \eta_1}{2} \norm{\x^0 - \x}^2 + \frac{\theta_1 \tau_1}{2} \norm{\y^0}^2  \nonumber \\
   & \quad + \inner{\y}{\theta_N \L(\x^N - \x^{N-1}) + \theta_1 \tau_1 (\y^N - \y^0)} + \sum_{k=1}^N \sum_{i=1}^m \theta_k e_i^k
 \end{align}
Then following the same argument as in the convex case. Define $\s^N$ and $\overline{\s}^N$ as before, we then can conclude $ \L \overline{\x}^N  = \overline{\s}^N $ and:
\begin{align}
  F(\overline{\x}^N) - F(\x^*) & \leqslant (\sum_{k=1}^N \theta_k)^{-1}
\left( \frac{\theta_1 \eta_1}{2} \norm{\x^0 - \x^*}^2 + \frac{\theta_1 \tau_1}{2} \norm{\y^0}^2  +\sum_{k=1}^N \sum_{i=1}^m \theta_k e_i^k\right ) \nonumber
\end{align}
Plug in definition of $\alpha_k = \frac{k}{k+1}, \theta_k = k+1, \eta_k = \frac{ku}{2}, \tau_k = \frac{4\norm{L}^2}{(k+1)u}$ and $e_i^k = \frac{\max(u \norm{\x^0 - \x^*}^2, \norm{L}^2 \norm{\y^0}^2/u)}{mNk}$ yields our convergence result.

\section{Complexity of DCGS: General Constraints}\label{complexity}
\subsection{Proof of Corollary \ref{complexity:smooth}}
\begin{proof}
From Theorem \ref{thrm:smooth}
, we know that to get an $\epsilon$-optimal solution, we need at most $N = \mathcal{O} \left(\frac{\norm{L}}{\epsilon} \max(\norm{\x^0 - \x^*}^2, \norm{\y^0}^2) \right)$ number of outer iterations in DCGS. Now we bound the number of calls to LO in the $k$-th outer iteration. Recall the CG procedure could be seemed as solving the subproblem in Line \ref{lo_subproblem} of Algorithm \ref{alg:dcgs} by Frank-Wolfe algorithm. From the well known result \citep{jaggi_revisiting_fw} we know that for solving a $l$-smooth function using Frank-Wolfe algorithm, with terminating wolfe-gap being $\epsilon$, the number of iterations could be bounded by $\frac{lD^2}{\epsilon}$. Observe that the objective function $\phi_i^k $ in the subproblem has smoothness $l+ \eta_k$, the total LO calls for each agent could be bounded by:
\begin{align}
 &\sum_{k=1}^N  \frac{(l+\eta_k)D^2}{e_i^k} \nonumber \\
   =  & \sum_{k=1}^N  \frac{(l + \norm{L}) D^2 mN}{\max(\norm{\x^0 - \x^*}^2, \norm{\y^0}^2)} \nonumber \\
   = & \frac{(l + \norm{L}) D^2 m N^2}{\max(\norm{\x^0 - \x^*}^2, \norm{\y^0}^2)} \nonumber \\
   = & \mathcal{O} \left( \frac{(l + \norm{L}) D^2 m \max(\norm{\x^0 - \x^*}^2, \norm{\y^0}^2)}{\epsilon^2} \right) \label{gen_lo_convex}
\end{align}
\end{proof}

\subsection{Proof of Corollary \ref{complexity:sc}}
\begin{proof}
From Theorem \ref{thrm:sc}, we know that to get an $\epsilon$-optimal solution, we need at most $N = \mathcal{O} \left(\sqrt{\frac{\max(u\norm{\x^0 - \x^*}^2, \frac{\norm{L}^2\norm{\y^0}^2}{u})}{\epsilon}}  \right)$ number of outer iterations in DCGS. Following the same argument as in
(\ref{gen_lo_convex}), we could be the LO complexity for each agent by:
\begin{align}
 &\sum_{k=1}^N  \frac{(l+\eta_k)D^2}{e_i^k} \nonumber \\
   =  & \sum_{k=1}^N  \frac{(l + ku) D^2 mNk}{\max(u\norm{\x^0 - \x^*}^2, \frac{\norm{L}^2\norm{\y^0}^2}{u})} \nonumber \\
   = & \sum_{k=1}^N  \frac{(l + ku) D^2 mNk}{\max(u\norm{\x^0 - \x^*}^2, \frac{\norm{L}^2\norm{\y^0}^2}{u})} \nonumber \\
   \leqslant & \mathcal{O} \left( \sum_{k=1}^N
      \frac{l D^2 m N k^2}{\max(u\norm{\x^0 - \x^*}^2, \frac{\norm{L}^2\norm{\y^0}^2}{u})} \right) \nonumber \\
      = & \mathcal{O} \left(
         \frac{l D^2 m N^4}{\max(u\norm{\x^0 - \x^*}^2, \frac{\norm{L}^2\norm{\y^0}^2}{u})} \right) \nonumber \\
  = & \mathcal{O} \left(
     \frac{l  m D^2\max(u\norm{\x^0 - \x^*}^2, \frac{\norm{L}^2\norm{\y^0}^2}{u})}{\epsilon^2} \right) \nonumber
\end{align}
\end{proof}

\section{Complexity of DCGS: Polyhedral Constraints}\label{complexity:poly}
\subsection{Proof of Corollary \ref{complexity:smooth_poly}}
\begin{proof}
From Theorem \ref{thrm:smooth} we know that to get an $\epsilon$-optimal solution, we need at most $N = \mathcal{O} \left(\frac{\norm{L}}{\epsilon} \max(\norm{\x^0 - \x^*}^2, \norm{\y^0}^2) \right)$ number of outer iteration in DCGS. Now we bound the number of calls to LO  in the  $k$-th outer iteration.
Our PCG procedure in DCGS could be seemed as pairwise Frank-Wolfe algorithm in \citep{julien_linear_convergence} applied to subproblem $\min_{x_i \in X} \phi_i^k(x_i)$. It has been shown that for a $u$-strongly convex and $l$-smooth function over a polyhedral set that has width $D$ and pyramidal width $W$, pairwise FW achieves a linear convergence rate of the wolfe gap. Specifically if we let $g^t$ denotes the wolfe gap at $t$-th iteration of pairwise FW algorithm, we have:
\begin{align}
  g^t \leqslant (1 - \frac{u}{4l} (\frac{W}{D})^2)^{\frac{t}{2}}\left(f(x^0) - f(x^*) \right)
\end{align}
Observe that at the $k$-th iteration of DCGS,
the subproblem in Line \ref{lo_subproblem} of DCGS has objective that is $l+\eta_k$ smooth and $\eta_k$ strongly convex, hence we could bound the LO of each agent by:
\begin{align}
& \mathcal{O} \left(\sum_{k=1}^N  \frac{(L + \eta_k)m D^2}{\eta_kW^2} \log(\frac{1}{e_i^k}) \right) \nonumber \\
= &  \mathcal{O} \left(\sum_{k=1}^N \frac{(L + \norm{L})m D^2}{\norm{L} W^2} \log(\frac{mN}{\norm{L} \max(\norm{\x^0 - \x^*}, \norm{\y^0}^2)}) \right) \nonumber \\
= & \tilde{\mathcal{O}} \left( \frac{(L + \norm{L})m D^2}{\norm{L} W^2} N \log(\frac{N m}{\norm{L}}) \right) \nonumber \\
= & \tilde{\mathcal{O}} \left(\frac{D^2}{W^2} \frac{\log m (l+\norm{L}) \max(\norm{\x^0 - \x^*}, \norm{\y^0}^2)}{\epsilon}  \right) \label{lo_sc_poly}
\end{align}
\end{proof}

\subsection{Proof of Corollary \ref{complexity:sc_poly}}
\begin{proof}
From Theorem \ref{thrm:sc}, we know that to get an $\epsilon$-optimal solution, we need at most $N = \mathcal{O} \left(\sqrt{\frac{\max(u\norm{\x^0 - \x^*}^2, \frac{\norm{L}^2\norm{\y^0}^2}{u})}{\epsilon}}  \right)$ number of outer iteration of DCGS. Following the same argument as in
(\ref{lo_sc_poly}), and note that the objective in subproblem $\min_{x_i \in X} \phi_i^k(x_i)$ is
$l+\eta_k$ smooth and $u+\eta_k$ strongly convex,
we could bound the LO complexity of each agent by:
\begin{align}
& \mathcal{O} \left(\sum_{k=1}^N \sum_{i=1}^m \frac{(L+\eta_k)D^2}{(u+\eta_k)W^2} \log(\frac{1}{e_i^k})\right) \nonumber \\
\leqslant & \mathcal{O} \left(\sum_{k=1}^N \frac{mLD^2}{uW^2} \log(\frac{mNk}{\max(u\norm{\x^0 - \x^*}^2, \frac{\norm{L}^2\norm{\y^0}^2}{u})})\right) \nonumber\\
= & \tilde{\mathcal{O}} \left(\frac{mLD^2}{uW^2} N\log(Nm) \right)\nonumber \\
= & \tilde{\mathcal{O}} \left( \frac{D^2}{W^2} \frac{l\log m}{\sqrt{u}} \frac{\norm{\x^0 -\x^*}}{\sqrt{\epsilon}} \right)
\end{align}
\end{proof}

\end{document}